\documentclass[12pt,diff1_eng]{article} 
\usepackage[english]{babel}
\usepackage{amssymb,amscd,amsmath}     
\makeatletter
\input{DIFF1_ENG.STY}
\makeatother
\usepackage{graphicx}
\usepackage{caption}
\usepackage[linktocpage=true, colorlinks=true, linkcolor=blue, citecolor=blue, urlcolor=blue]{hyperref}

\begin{document}

\large  


\begin{center}
\Large{\bf New representations of the Hu--Meyer formulas 
and series expansion of iterated
Stratonovich stochastic integrals with
respect to components of a multidimensional Wiener process}
\end{center}

$$
{}
$$

\centerline{\large {Dmitriy F. Kuznetsov}}

\begin{center}
Peter the Great Saint-Petersburg Polytechnic University\\
Polytechnicheskaya ul., 29,\\
195251, Saint-Petersburg, Russia\\
\hspace{-2mm}e-mail: sde\_\hspace{0.5mm}kuznetsov@inbox.ru
\end{center}

$$
{}
$$

\noindent
{\bf Abstract.} The article is devoted to the systematic derivation of new 
representations of the Hu--Meyer formulas.
The formula expressing a multiple Wiener stochastic integral through the sum of multiple 
Stratonovich stochastic integrals and the formula expressing a multiple Stratonovich 
stochastic integral through the sum of multiple Wiener stochastic integrals are derived
for the case of a 
multidimensional Wiener process.
At that several different definitions of the multiple Stratonovich stochastic
integral and several variants of sufficient conditions for the validity of the 
Hu--Meyer formulas are used. In particular, the proof method proposed by the author in 2006 is applied
to obtain
Hu--Meyers formulas based on generalized multiple Fourier series
for the case of a multidimensional Wiener process.
Of great importance for the numerical solution 
of It\^{o} stochastic differential equations
is the verification of sufficient conditions for the applicability of the 
Hu--Meyer formula (based on generalized multiple Fourier series)
for the case of iterated Stratonovich stochastic integrals 
with respect to components of a multidimensional Wiener process.
In the author's previous works, the indicated conditions were verified 
for iterated Stratonovich stochastic integrals of multiplicities 1 to 6  
(the case of an arbitrary basis in the Hilbert space) 
and for iterated Stratonovich stochastic integrals of multiplicities 7 and 8
(the case of two special bases in 
Hilbert space (the trigonometric Fourier basis and the basis of Legendre polynomials)).
The theorems on the expansion of iterated Stratonovich stochastic integrals 
of arbitrary multiplicity are proved.
These theorems prove Hypotheses 2.1, 2.2, and 2.5 from the monograph \cite{2018a} (Sect.~2.5, 2.28).

\renewcommand{\baselinestretch}{1.3}

\vspace{10mm}

{\normalsize 

\linespread{1.0}

\tableofcontents

\vspace{8mm}

\linespread{1.0}

}





\section{Introduction}

In the theory of stochastic integrals, 
a special place is occupied by the so-called Hu--Meyer formulas, 
which express a multiple Wiener stochastic integral through the sum of multiple 
Stratonovich stochastic integrals and, conversely, express a multiple Stratonovich 
stochastic integral through the sum of multiple Wiener stochastic integrals.
Usually, the indicated multiple stochastic integrals in the Hu--Meier formulas are 
considered with respect to a scalar Wiener or Poisson processes \cite{ito1951}-\cite{hida}.
It should be noted that a version of the Hu--Meyer formula for multiple 
stochastic integrals with respect to L{\'e}vy processes is known \cite{farre}.

Of great interest (in the context of numerical integration of It\^{o} stochastic 
differential equations (SDEs)) is the Hu--Meyer formula, which express the 
iterated It\^{o} stochastic integral with respect to the components of a 
multidimensional Wiener process through the sum of multiple Stratonovich stochastic integrals.
A version of such a formula based on generalized multiple Fourier series was 
first obtained in 2006 \cite{2006} (also see \cite{2009}-\cite{arxiv-9}).
The inverse version of the above formula, which expresses the multiple 
(or iterated) Stratonovich stochastic integral in terms of multiple Wiener 
stochastic integrals, can be found in \cite{Rybakov30001}.

Rest of the work is organized as follows. 
In Sect.~2.1 we give  
a brief review of some old results on the expansion
of iterated Stratonovich and It\^{o} stochastic integrals obtained by the author.
Sect.~3.1 and 3.2 are devoted to the systematic derivation of new 
representations for both of the indicated variants of the Hu--Meyer formula: 
the formula expressing a multiple Wiener stochastic integral through the sum of multiple 
Stratonovich stochastic integrals, and the formula expressing a multiple Stratonovich 
stochastic integral through the sum of multiple Wiener stochastic integrals (the case of a 
multidimensional Wiener process).
Several different definitions of 
the multiple Stratonovich stochastic integral and several variants of sufficient conditions 
for the validity of the Hu--Meyer formulas are used. 
At that, we use the proof method from our works \cite{2006}, 
\cite{arxiv-1}, \cite{diffjournal} to obtain
Hu--Meyers formulas based on genealized multiple Fourier series
for the case of a multidimensional Wiener process.

It should be noted that of great importance for the numerical integration 
of It\^{o} SDEs is the verification of sufficient conditions for the applicability of the 
Hu--Meyer formula (based on generalized multiple Fourier series)
for the case of iterated Stratonovich stochastic integrals 
with respect to components of a multidimensional Wiener process.
In the author's previous works \cite{new-art-1xxys1}-\cite{2018a}
the indicated conditions were verified 
for iterated Stratonovich stochastic integrals of multiplicities 1 to 6  
(the case of an arbitrary basis in the Hilbert space)
and for iterated Stratonovich stochastic integrals of multiplicities 7 and 8
(the case of two special bases in 
Hilbert space (the trigonometric Fourier basis and the basis of Legendre polynomials)).

The theorems on the expansion of iterated Stratonovich stochastic integrals 
of arbitrary multiplicity are proved (Sect.~3.3).

\section{Preliminary Results}

\subsection{A Review of Some Old Results on the Expansion
of Iterated Stratonovich and It\^{o} Stochastic Integrals}

In this section, we will give  
a brief review of some old results on the expansion
of iterated Stratonovich and It\^{o} stochastic integrals
obtained by the author (the case of an arbitrary complete orthonormal
system (CONS) in 
$L_2([t, T])$).

Let $(\Omega,{\rm F},{\sf P})$ be a complete probability space, let 
$\{{\rm F}_t, t\in[0,T]\}$ be a nondecreasing 
right-continous family of $\sigma$-algebras of ${\rm F},$
and let ${\bf w}_t=({\bf w}_{t}^{(1)}\ldots {\bf w}_{t}^{(m)})^{\sf T}$ be 
a standard $m$-dimensional Wiener 
stochastic process with independent components which is
${\rm F}_t$-measurable for any $t\in[0, T].$ 
Consider the following
iterated Stratonovich stochastic integrals
\begin{equation}
\label{str}
J^{*}[\psi^{(k)}]_{T,t}^{(i_1 \ldots i_k)}=
{\int\limits_t^{*}}^T
\psi_k(t_k) \ldots 
{\int\limits_t^{*}}^{t_2}
\psi_1(t_1) d{\bf w}_{t_1}^{(i_1)}\ldots
d{\bf w}_{t_k}^{(i_k)},
\end{equation}
where $\psi_1(\tau),\ldots,\psi_k(\tau): [t, T]\to\mathbb{R},$
$i_1,\ldots,i_k = 0,$ $1,$ $\ldots,$ $m,$
${\bf w}_{\tau}^{(0)}=\tau$.

Note that the importance of the problem of numerical integration of 
It\^{o} SDEs is explained by a wide range of
their applications \cite{KlPl2}-\cite{KPS}.               
One of the effective approaches to the above
problem is based on 
the so-called Taylor--Stratonovich expansions \cite{KlPl2} (also see \cite{2018a}). 
The stochastic integrals (\ref{str}) are multipliers for
terms of the mentioned expansions
that explains the relevance of the topic of the work.

Throughout this article we will use the definition of the iterated Stratonovich stochastic
integral from \cite{KlPl2} (also see \cite{2018a}, Sect.~2.1.1).

Let $\psi_1(\tau),\ldots,\psi_k(\tau)\in L_2([t, T]).$
Consider the Fourier coefficient
\begin{equation}
\label{after3000}
C_{j_k \ldots j_1}=\int\limits_t^T\psi_k(t_k)\phi_{j_k}(t_k)\ldots
\int\limits_t^{t_2}
\psi_1(t_1)\phi_{j_1}(t_1)
dt_1\ldots dt_k
\end{equation}

\noindent
corresponding to the factorized Volterra--type kernel of the form
\begin{equation}
\label{july7000}
K(t_1,\ldots,t_k)=
\left\{\begin{matrix}
\psi_1(t_1)\ldots \psi_k(t_k),\ &t_1<\ldots<t_k\cr\cr
0,\ &\hbox{\rm otherwise}
\end{matrix}
\right.,
\end{equation}

\noindent
where $(t_1,\ldots,t_k)\in [t, T]^k$ $(k\ge 2)$ and
$\{\phi_j(x)\}_{j=0}^{\infty}$ is
an arbitrary CONS in 
$L_2([t, T])$.

{\bf Theorem~1}\ \cite{2018a} (Sect.~2.1.4).\ {\it Suppose that 
$\{\phi_j(x)\}_{j=0}^{\infty}$ is an arbitrary CONS 
in $L_2([t, T]).$
Moreover$,$ $\psi_1(\tau), \psi_2(\tau)$ are continuous 
functions on $[t, T].$ 
Then
\begin{equation}
\label{trace20}
J^{*}[\psi^{(2)}]_{T,t}^{(i_1 i_2)}=\hbox{\vtop{\offinterlineskip\halign{
\hfil#\hfil\cr
{\rm l.i.m.}\cr
$\stackrel{}{{}_{p_1,p_2\to \infty}}$\cr
}} }\sum_{j_1=0}^{p_1}\sum_{j_2=0}^{p_2}
C_{j_2j_1}\zeta_{j_1}^{(i_1)}\zeta_{j_2}^{(i_2)},
\end{equation}
where $J^{*}[\psi^{(2)}]_{T,t}^{(i_1 i_2)}$ is defined by {\rm (\ref{str}),}
$i_1, i_2=0, 1,\ldots,m,$ 
${\rm l.i.m.}$ is a limit in the mean-square sense$,$
$C_{j_2 j_1}$ has the form {\rm (\ref{after3000}),}
$$
\zeta_{j}^{(i)}=
\int\limits_t^T \phi_{j}(\tau) d{\bf w}_{\tau}^{(i)}
$$ 
are independent standard Gaussian random variables for various 
$i$ or $j$ {\rm (}in the case when $i\ne 0${\rm ),}
and ${\bf w}_{\tau}^{(0)}=\tau.$}

{\bf Theorem~2}\ \cite{2024xxx}, \cite{2018a} (Sect.~2.24, 2.33).\ {\it Suppose that
$\{\phi_j(x)\}_{j=0}^{\infty}$ is an arbitrary CONS 
in $L_2([t,T])$ and $\psi_i(\tau)=(\tau-t)^{l_i}$
$(l_i=0,1,2,\ldots,\ i=1,\ldots,4)$.
Then
$$
J^{*}[\psi^{(k)}]_{T,t}^{(i_1 \ldots i_k)}
=
\hbox{\vtop{\offinterlineskip\halign{
\hfil#\hfil\cr
{\rm l.i.m.}\cr
$\stackrel{}{{}_{p\to \infty}}$\cr
}} }\sum_{j_1,\ldots,j_k=0}^{p}
C_{j_k \ldots j_1}\zeta_{j_1}^{(i_1)}\ldots \zeta_{j_k}^{(i_k)}\ \ \ (k=3,\ 4),
$$
where $J^{*}[\psi^{(k)}]_{T,t}^{(i_1 \ldots i_k)}$ $(k=3,\ 4)$
is defined by {\rm (\ref{str}),} 
$i_1,\ldots,i_4=0,1,\ldots,m,$ $C_{j_k \ldots j_1}$ has the form {\rm (\ref{after3000});} 
another notations as in Theorem~{\rm 1.}}

{\bf Theorem~3}\ \cite{2024xxx}, \cite{2018a} (Sect.~2.32, 2.34).\ {\it Suppose that
$\{\phi_j(x)\}_{j=0}^{\infty}$ is an arbitrary CONS 
in $L_2([t,T])$ and $\psi_1(\tau),\ldots,\psi_6(\tau)\equiv 1.$
Then
$$
J^{*}[\psi^{(k)}]_{T,t}^{(i_1\ldots i_k)}=
\hbox{\vtop{\offinterlineskip\halign{
\hfil#\hfil\cr
{\rm l.i.m.}\cr
$\stackrel{}{{}_{p\to \infty}}$\cr
}} }
\sum\limits_{j_1,\ldots, j_k=0}^{p}
C_{j_k \ldots j_1}\zeta_{j_1}^{(i_1)}\ldots \zeta_{j_k}^{(i_k)}\ \ \ (k=5,\ 6),
$$
where $J^{*}[\psi^{(k)}]_{T,t}^{(i_1\ldots i_k)}$ $(k=5,\ 6)$
is defined by {\rm (\ref{str}),}
$i_1,\ldots,i_6=0, 1,\ldots,m,$ $C_{j_k \ldots j_1}$ has the form {\rm (\ref{after3000});}
another notations as in Theorem~{\rm 1.}}

Consider the unordered
set $\{1, 2, \ldots, k\}$ 
and separate it into two parts:
the first part consists of $r$ unordered 
pairs (sequence order of these pairs is also unimportant) and the 
second one consists of the 
remaining $k-2r$ numbers.
So, we have 
\begin{equation}
\label{leto5007after}
(\{
\underbrace{\{g_1, g_2\}, \ldots, 
\{g_{2r-1}, g_{2r}\}}_{\small{\hbox{part 1}}}
\},
\{\underbrace{q_1, \ldots, q_{k-2r}}_{\small{\hbox{part 2}}}
\}),
\end{equation}

\noindent
where 
$\{g_1, g_2, \ldots, 
g_{2r-1}, g_{2r}, q_1, \ldots, q_{k-2r}\}=\{1, 2, \ldots, k\},$
braces   
mean an unordered 
set, and pa\-ren\-the\-ses mean an ordered set.

Denote
$$
C_{j_k \ldots j_{l+1}j_l j_{l-1} j_{l-2} \ldots j_1}\biggl|_{(j_l j_{l-1})
\curvearrowright (\cdot) }\biggr.
\stackrel{\sf def}{=}
$$

\vspace{1mm}
$$
\stackrel{\sf def}{=}\int\limits_t^T\psi_k(t_k)\phi_{j_k}(t_k)\ldots
\int\limits_t^{t_{l+2}}\psi_{l+1}(t_{l+1})\phi_{j_{l+1}}(t_{l+1})
\int\limits_t^{t_{l+1}}\psi_{l}(t_{l})\psi_{l-1}(t_{l})\times
$$

\vspace{-3mm}
\begin{equation}
\label{after900}
\times
\int\limits_t^{t_{l}}\psi_{l-2}(t_{l-2})\phi_{j_{l-2}}(t_{l-2})\ldots 
\int\limits_t^{t_2}
\psi_1(t_1)\phi_{j_1}(t_1)
dt_1\ldots dt_{l-2}dt_{l}t_{l+1}\ldots dt_k,
\end{equation}

\vspace{1mm}
\noindent
where $\{l, l-1\}$ is one of the pairs
$\{g_1, g_2\}, \ldots,\{g_{2r-1}, g_{2r}\}$ $($see {\rm (\ref{leto5007after})}$)$.

Further, we will use the quantity
$$
\prod\limits_{l=1}^r {\bf 1}_{\{g_{2l}=g_{2l-1}+1\}}
C_{j_k \ldots j_1}\biggl|_{(j_{g_2} j_{g_1})\curvearrowright (\cdot)
\ldots (j_{g_{2r}} j_{g_{2r-1}})\curvearrowright (\cdot),
j_{g_{{}_{1}}}=~j_{g_{{}_{2}}},\ldots, j_{g_{{}_{2r-1}}}=~j_{g_{{}_{2r}}}
}\biggr.
$$

\vspace{1mm}
\noindent
which is defined iteratively using the equality (\ref{after900}),
where ${\bf 1}_A$ is the indicator of the set $A$.

Introduce the following notations
$$
J[\psi^{(k)}]_{T,t}^{(i_1\ldots i_k)[s_l,\ldots,s_1]}\ \stackrel{\rm def}{=}\ 
\prod_{q=1}^l {\bf 1}_{\{i_{s_q}=
i_{s_{q}+1}\ne 0\}}\ \times
$$
$$
\times
\int\limits_t^T\psi_k(t_k)\ldots \int\limits_t^{
t_{s_l+3}}\psi_{s_l+2}(t_{s_l+2})
\int\limits_t^{t_{s_l+2}}\psi_{s_l}(t_{s_l+1})
\psi_{s_l+1}(t_{s_l+1})\times
$$
$$
\times
\int\limits_t^{t_{s_l+1}}\psi_{s_l-1}(t_{s_l-1})
\ldots
\int\limits_t^{t_{s_1+3}}\psi_{s_1+2}(t_{s_1+2})
\int\limits_t^{t_{s_1+2}}\psi_{s_1}(t_{s_1+1})
\psi_{s_1+1}(t_{s_1+1})\times
$$
$$
\times
\int\limits_t^{t_{s_1+1}}\psi_{s_1-1}(t_{s_1-1})
\ldots \int\limits_t^{t_2}\psi_1(t_1)
d{\bf w}_{t_1}^{(i_1)}\ldots d{\bf w}_{t_{s_1-1}}^{(i_{s_1-1})}
dt_{s_1+1}d{\bf w}_{t_{s_1+2}}^{(i_{s_1+2})}\ldots
$$
\begin{equation}
\label{30.1}
\ldots
d{\bf w}_{t_{s_l-1}}^{(i_{s_l-1})}
dt_{s_l+1}d{\bf w}_{t_{s_l+2}}^{(i_{s_l+2})}\ldots d{\bf w}_{t_k}^{(i_k)},
\end{equation}

\vspace{4mm}
\noindent
where $\psi_1(\tau),\ldots,\psi_k(\tau)\in L_2([t, T]),$
the iterated integral on the right-hand side of
(\ref{30.1}) is 
the iterated It\^{o} stochastic integral,
$(s_l,\ldots,s_1)\in{\rm A}_{k,l},$

\vspace{-5mm}
\begin{equation}
\label{30.5550001}
{\rm A}_{k,l}=\bigl\{(s_l,\ldots,s_1):\
s_l>s_{l-1}+1,\ldots,s_2>s_1+1;\ s_l,\ldots,s_1=1,\ldots,k-1\bigr\},
\end{equation}

\noindent
$l=1,2,\ldots,\left[k/2\right],$\ 
$i_1,\ldots,i_k=0, 1,\ldots,m,$\ 
$[x]$ is an
integer part of a real number $x,$\ 
${\bf 1}_A$ is the indicator of the set $A$.

For $l=0,$ we understand (\ref{30.1}) as the following
iterated It\^{o} stochastic integral
\begin{equation}
\label{ito}
J[\psi^{(k)}]_{T,t}^{(i_1 \ldots i_k)}=\int\limits_t^T\psi_k(t_k) \ldots 
\int\limits_t^{t_{2}}
\psi_1(t_1) d{\bf w}_{t_1}^{(i_1)}\ldots
d{\bf w}_{t_k}^{(i_k)},
\end{equation}
where 
$i_1,\ldots,i_k = 0,1,\ldots,m$ and
${\bf w}_{\tau}^{(0)}=\tau$.

Let us denote
\begin{equation}
\label{dsds9}
J[\psi^{(k)}]_{T,t}^{(i_1\ldots i_k)}+
\sum_{r=1}^{\left[k/2\right]}\frac{1}{2^r}
\sum_{(s_r,\ldots,s_1)\in {\rm A}_{k,r}}
J[\psi^{(k)}]_{T,t}^{(i_1\ldots i_k)[s_r,\ldots,s_1]}
\stackrel{\sf def}{=}
\bar J^{*}[\psi^{(k)}]_{T,t}^{(i_1\ldots i_k)},
\end{equation}

\noindent
where $\psi_1(\tau),\ldots,\psi_k(\tau)\in L_2([t, T])$,
$J[\psi^{(k)}]_{T,t}^{(i_1\ldots i_k)}$ is the iterated It\^{o} stochastic
integral (\ref{ito}),
$\sum\limits_{\emptyset}$ is supposed to be equal to zero.

{\bf Theorem~4}\ \cite{2018a} (Sect.~2.22), \cite{arxiv-4}, \cite{arxiv-5}.\ {\it Assume that
the CONS $\{\phi_j(x)\}_{j=0}^{\infty}$
in $L_2([t, T])$ and
$\psi_1(\tau),\ldots, \psi_k(\tau)\in L_2([t, T])$
are such that 
$$
\lim\limits_{p\to\infty}
\sum\limits_{j_{q_1},\ldots,j_{q_{k-2r}}=0}^p
\Biggl(\sum\limits_{j_{g_1}, j_{g_3},\ldots ,j_{g_{2r-1}}=0}^p
C_{j_k\ldots j_1}\biggl|_{j_{g_1}=j_{g_2},\ldots, j_{g_{2r-1}}=j_{g_{2r}}}-\Biggr.
$$

\vspace{-4mm}
\begin{equation}
\label{july700000}
\Biggl.-\frac{1}{2^r} \prod\limits_{l=1}^r {\bf 1}_{\{g_{2l}=g_{2l-1}+1\}}
C_{j_k \ldots j_1}\biggl|_{(j_{g_2} j_{g_1})\curvearrowright (\cdot)
\ldots (j_{g_{2r}} j_{g_{2r-1}})\curvearrowright (\cdot),
j_{g_{{}_{1}}}=~j_{g_{{}_{2}}},\ldots, j_{g_{{}_{2r-1}}}=~j_{g_{{}_{2r}}}
}\biggr.\Biggr)^2=0
\end{equation}

\vspace{1mm}
\noindent
for all $r=1, 2,\ldots,[k/2]$ 
and for all possible $g_1,g_2,\ldots,g_{2r-1},g_{2r}$ {\rm (}see {\rm (\ref{leto5007after}))}.
Then$,$ for the sum $\bar J^{*}[\psi^{(k)}]_{T,t}^{(i_1\ldots i_k)}$
of iterated It\^{o} stochastic integrals 
defined by {\rm (\ref{dsds9})}
the following 
expansion 
$$
\bar J^{*}[\psi^{(k)}]_{T,t}^{(i_1\ldots i_k)}=
\hbox{\vtop{\offinterlineskip\halign{
\hfil#\hfil\cr
{\rm l.i.m.}\cr
$\stackrel{}{{}_{p\to \infty}}$\cr
}} }
\sum_{j_1,\ldots,j_k=0}^{p}
C_{j_k \ldots j_1}\prod\limits_{l=1}^k \zeta_{j_l}^{(i_l)}
$$

\noindent
that converges in the mean-square sense is valid$,$ where 
$C_{j_k \ldots j_1}$ has the form {\rm (\ref{after3000}),}
${\rm l.i.m.}$ is a limit in the mean-square sense$,$
$i_1, \ldots, i_k=0, 1,\ldots,m,$
$$
\zeta_{j}^{(i)}=
\int\limits_t^T \phi_{j}(\tau) d{\bf w}_{\tau}^{(i)}
$$ 
are independent standard Gaussian random variables for various 
$i$ or $j$ {\rm (}in the case when $i\ne 0${\rm ),}
and 
${\bf w}_{\tau}^{(0)}=\tau.$}

Let us formulate the statements on connection 
between
iterated 
Stra\-to\-no\-vich and It\^{o} stochastic integrals 
(\ref{str}) and (\ref{ito})
of arbitrary multiplicity $k$ $(k\in{\bf N})$.

{\bf Theorem 5} \cite{300a} (1997) (also see \cite{2018a} (Sect.~2.4.1, Theorem~2.12)).\
{\it Suppose that
$\psi_1(\tau),$ $\ldots,\psi_k(\tau)$ are continuous
functions on $[t, T]$.
Then
\begin{equation}
\label{30.4}
J^{*}[\psi^{(k)}]_{T,t}^{(i_1\ldots i_k)}=J[\psi^{(k)}]_{T,t}^{(i_1\ldots i_k)}+
\sum_{r=1}^{\left[k/2\right]}\frac{1}{2^r}
\sum_{(s_r,\ldots,s_1)\in {\rm A}_{k,r}}
J[\psi^{(k)}]_{T,t}^{(i_1\ldots i_k)[s_r,\ldots,s_1]}\ \ \ \hbox{{\rm w.\ p.\ 1}},
\end{equation}
where $i_1,\ldots,i_k=0, 1,\ldots,m,$ $\sum\limits_{\emptyset}$ is supposed to be equal to zero,
w.~p.~$1$ {\rm (}here and further{\rm )} means with probability $1$.}

{\bf Theorem 6} \cite{2018a} (Sect.~4.6, Proposition~4.3).\
{\it Suppose that
$\psi_1(\tau),$ $\ldots,\psi_k(\tau)$ are continuous
functions on $[t, T]$.
Then
\begin{equation}
\label{2024str11}
J[\psi^{(k)}]_{T,t}^{(i_1\ldots i_k)}=J^{*}[\psi^{(k)}]_{T,t}^{(i_1\ldots i_k)}+
\sum_{r=1}^{\left[k/2\right]}\frac{(-1)^r}{2^r}
\sum_{(s_r,\ldots,s_1)\in {\rm A}_{k,r}}
J^{*}[\psi^{(k)}]_{T,t}^{(i_1\ldots i_k)[s_r,\ldots,s_1]}\ \ \hbox{{\rm w.\ p.\ 1,}}
\end{equation}
where $J^{*}[\psi^{(k)}]_{T,t}^{(i_1\ldots i_k)[s_r,\ldots,s_1]}$
is defined similarly to $J[\psi^{(k)}]_{T,t}^{(i_1\ldots i_k)[s_r,\ldots,s_1]},$
only on the right-hand side of {\rm (\ref{30.1})} all It\^{o} integrals
are replaced by Stratonovich integrals$;$
another notations as in Theorem~{\rm 5}.}

Note that the condition of continuity of the functions
$\psi_1(\tau),\ldots, \psi_k(\tau)$ 
is related to the definition  
of the Stratonovich stochastic integral that we use
\cite{KlPl2} (also see \cite{2018a}, Sect.~2.1.1).

Combining Theorems~4 and 5, we obtain the following result.

{\bf Theorem 7}\ \cite{2018a} (Sect.~2.22), \cite{arxiv-4}, \cite{arxiv-5}.\
{\it Suppose that
the CONS $\{\phi_j(x)\}_{j=0}^{\infty}$
in $L_2([t, T])$ and
continuous functions $\psi_1(\tau),\ldots, \psi_k(\tau)$ on $[t, T]$
are such that the condition {\rm (\ref{july700000})} is satisfied.
Then
\begin{equation}
\label{asasas}
J^{*}[\psi^{(k)}]_{T,t}^{(i_1\ldots i_k)}=
\hbox{\vtop{\offinterlineskip\halign{
\hfil#\hfil\cr
{\rm l.i.m.}\cr
$\stackrel{}{{}_{p\to \infty}}$\cr
}} }
\sum_{j_1,\ldots,j_k=0}^{p}
C_{j_k \ldots j_1}\prod\limits_{l=1}^k \zeta_{j_l}^{(i_l)},
\end{equation}
where $J^{*}[\psi^{(k)}]_{T,t}^{(i_1\ldots i_k)}$
is the iterated Stratonovich stochastic integral {\rm (\ref{str});}
another notations as in Theorem~{\rm 4}.}

Suppose that $k>2r.$  Consider the following variant 
of sufficient condition under which the condition (\ref{july700000}) is satisfied 
(see \cite{2018a}, Sect.~2.29, 2.30)

\vspace{-3mm}
$$
\exists\ \ \ \lim\limits_{p,q\to\infty}
\sum\limits_{j_{q_1},\ldots,j_{q_{k-2r}}=0}^q
\Biggl(\sum\limits_{j_{g_1}, j_{g_3},\ldots ,j_{g_{2r-1}}=0}^p
C_{j_k\ldots j_1}\biggl|_{j_{g_1}=j_{g_2},\ldots, j_{g_{2r-1}}=j_{g_{2r}}}-\Biggr.
$$

\vspace{-2mm}
\begin{equation}
\label{july700001xyz}
\Biggl.-\frac{1}{2^r} \prod\limits_{l=1}^r {\bf 1}_{\{g_{2l}=g_{2l-1}+1\}}
C_{j_k \ldots j_1}\biggl|_{(j_{g_2} j_{g_1})\curvearrowright (\cdot)
\ldots (j_{g_{2r}} j_{g_{2r-1}})\curvearrowright (\cdot),
j_{g_{{}_{1}}}=~j_{g_{{}_{2}}},\ldots, j_{g_{{}_{2r-1}}}=~j_{g_{{}_{2r}}}
}\biggr.\Biggr)^2<\infty
\end{equation}

\vspace{1mm}
\noindent 
for all $r=1,2,\ldots,[k/2],$
where notations are the same as in 
(\ref{july700000}).
Recall that the case
$k=2r$ of (\ref{july700000}) is proved in \cite{rybakov5000}-\cite{rybakov5001} 
(also see \cite{2024xxx}, Sect.~3.2 or \cite{2018a}, Sect.~2.27.4).

{\bf Theorem 8}\ \cite{2018a} (Sect.~1.11).\
{\it Suppose that 
$\psi_1(\tau),\ldots,\psi_k(\tau)\in L_2([t, T])$ and
$\{\phi_j(x)\}_{j=0}^{\infty}$ is an arbitrary CONS 
in $L_2([t,T]).$
Then 
\begin{equation}
\label{febr5000}
J[\psi^{(k)}]_{T,t}^{(i_1\ldots i_k)}=
\hbox{\vtop{\offinterlineskip\halign{
\hfil#\hfil\cr
{\rm l.i.m.}\cr
$\stackrel{}{{}_{p_1,\ldots,p_k\to \infty}}$\cr
}} }
\sum\limits_{j_1=0}^{p_1}\ldots
\sum\limits_{j_k=0}^{p_k}
C_{j_k\ldots j_1}J'[\phi_{j_1}\ldots \phi_{j_k}]_{T,t}^{(i_1\ldots i_k)},
\end{equation}

\noindent
where $i_1,\ldots,i_k=0,1,\ldots,m,$ 
$J[\psi^{(k)}]_{T,t}^{(i_1\ldots i_k)}$ is the iterated It\^{o}
stochastic integral {\rm (\ref{ito}),}
$J'[\phi_{j_1}\ldots \phi_{j_k}]_{T,t}^{(i_1\ldots i_k)}$ is the 
multiple Wiener stochastic integral
defined as in {\rm \cite{ito1951}} but for the case of a multidimensional Wiener
process {\rm (}see {\rm \cite{2018a},} Sect.~{\rm 1.11} for details{\rm ),}
$C_{j_k \ldots j_1}$ is the Fourier coefficient 
{\rm (\ref{after3000});} another notations as in Theorem~{\rm 1}.}

\section{Main Results}

\subsection{The Connection of Condition (\ref{july700000}) with the Concept
of Limiting Traces from the Work of G.W. Johnson and G. Kallianpur \cite{bugh3}}

Assume that $\{\phi_j(x)\}_{j=0}^{\infty}$
is an arbitrary CONS of functions
in $L_2([t, T])$ and
$\psi_1(\tau),\ldots, \psi_k(\tau)\in L_2([t, T]).$

By analogy with Definition~2.5 in \cite{bugh2}, we define the 
$r$th limiting trace of the function $K(t_1,\ldots,t_k)\in L_2([t, T]^k)$
of type (\ref{july7000}) by the following expression
$$
T^{k-2r}_{g_1,g_2,\ldots,g_{2r-1}, g_{2r}}
K (t_{q_1},\ldots,t_{q_{k-2r}})
\stackrel{\sf def}{=}
$$
\begin{equation}
\label{2025may22}
\stackrel{\sf def}{=}
\lim\limits_{p\to\infty}T^{k-2r,\hspace{0.2mm}p}_{g_1,g_2,\ldots,g_{2r-1}, g_{2r}}
K(t_{q_1},\ldots,t_{q_{k-2r}})
\end{equation}

\vspace{1mm}
\noindent
in $L_2([t, T]^{k-2r})$ (here and further $L_2([t, T]^{0})\stackrel{\sf def}{=}{\bf R}$), where

\vspace{-2mm}
$$
T^{k-2r,\hspace{0.2mm}p}_{g_1,g_2,\ldots,g_{2r-1}, g_{2r}}
K(t_{q_1},\ldots,t_{q_{k-2r}})
=
$$

\vspace{-4mm}
$$
=
\sum\limits_{j_{q_1},\ldots,j_{q_{k-2r}}=0}^p
\sum\limits_{j_{g_1}, j_{g_3},\ldots ,j_{g_{2r-1}}=0}^p
C_{j_k\ldots j_1}\biggl|_{j_{g_1}=j_{g_2},\ldots, j_{g_{2r-1}}=j_{g_{2r}}}\times
$$

\vspace{-1mm}
$$
\times\phi_{j_{q_1}}(t_{q_1})\ldots \phi_{j_{q_{k-2r}}}(t_{q_{k-2r}}),
$$

\vspace{3mm}
\noindent
where
$r=1,2,\ldots,[k/2]$ and 
$\{g_1,g_2,\ldots,g_{2r-1},g_{2r},
q_1,\ldots,q_{k-2r}\}=\{1,2,\ldots,k\}$ 
(see (\ref{leto5007after})),
$C_{j_k \ldots j_1}$ is the Fourier coefficient 
(\ref{after3000}).         
In addition we write
$T^{k-2r}_{g_1,g_2,\ldots,g_{2r-1}, g_{2r}}
K(t_{q_1},\ldots,t_{q_{k-2r}})
=K(t_1,\ldots,t_k)$ 
for $r=0$.

Note that in \cite{bugh3} the Wiener process is scalar,
while in this article the Wiener process is a multidimensional
process with independent components.
One of the main results of work \cite{bugh3} (Theorem~5.1)
is obtained under the condition 
of existence of limiting traces (see Definition~2.5 in \cite{bugh2}).

Further, we will show that the condition (\ref{july700000})
is a necessary and sufficient condition 
for the existence of limiting traces (\ref{2025may22})
(for all $r=1,2,\ldots,[k/2]$ and for all possible $g_1,g_2,\ldots,g_{2r-1},g_{2r}$ 
(see (\ref{leto5007after}))) for the case of a 
multidimensional Wiener process and $K(t_1,\ldots,t_k)$
defined by (\ref{july7000}).

Here it is also appropriate to recall the formula (see the proof in \cite{2018a}, Sect.~2.18)
$$
\hbox{\vtop{\offinterlineskip\halign{
\hfil#\hfil\cr
{\rm l.i.m.}\cr
$\stackrel{}{{}_{p_1,\ldots,p_k\to \infty}}$\cr
}} }\sum_{j_1=0}^{p_1}\ldots\sum_{j_k=0}^{p_k}
C_{j_k\ldots j_1}
\zeta_{j_1}^{(i_1)}\ldots \zeta_{j_k}^{(i_k)}
=J[\psi^{(k)}]_{T,t}^{(i_1\ldots i_k)}
+
$$

$$
+
\sum\limits_{r=1}^{[k/2]}
\sum_{\stackrel{(\{\{g_1, g_2\}, \ldots, 
\{g_{2r-1}, g_{2r}\}\}, \{q_1, \ldots, q_{k-2r}\})}
{{}_{\{g_1, g_2, \ldots, 
g_{2r-1}, g_{2r}, q_1, \ldots, q_{k-2r}\}=\{1, 2, \ldots, k\}}}}
\prod\limits_{s=1}^r
{\bf 1}_{\{i_{g_{{}_{2s-1}}}=~i_{g_{{}_{2s}}}\ne 0\}}\times
$$

\vspace{2mm}
\begin{equation}
\label{2025may23}
\times \hbox{\vtop{\offinterlineskip\halign{
\hfil#\hfil\cr
{\rm l.i.m.}\cr
$\stackrel{}{{}_{p_1,\ldots,p_k\to \infty}}$\cr
}} }\sum_{j_1=0}^{p_1}\ldots\sum_{j_k=0}^{p_k}
C_{j_k\ldots j_1}
\prod\limits_{s=1}^r{\bf 1}_{\{j_{g_{{}_{2s-1}}}=~j_{g_{{}_{2s}}}\}}
J'[\phi_{j_{q_1}}\ldots \phi_{j_{q_{k-2r}}}]_{T,t}^{(i_{q_1}\ldots i_{q_{k-2r}})}
\end{equation}

\vspace{1mm}
\noindent
w.~p.~1, where 
$J'[\phi_{j_{q_1}}\ldots \phi_{j_{q_{k-2r}}}]_{T,t}^{(i_{q_1}\ldots i_{q_{k-2r}})}$ 
is the 
multiple Wiener stochastic integral
defined as in \cite{ito1951} but for the case of a multidimensional Wiener
process (see \cite{2018a}, Sect.~1.11 for details),
$J'[\phi_{j_{q_1}}\ldots \phi_{j_{q_{k-2r}}}]_{T,t}^{(i_{q_1}\ldots i_{q_{k-2r}})}
\stackrel{\sf def}{=}1$ for $k=2r,$
$J[\psi^{(k)}]_{T,t}^{(i_1\ldots i_k)}$ is the iterated It\^{o} stochastic
integral (\ref{ito}), $C_{j_k \ldots j_1}$ is the Fourier coefficient 
(\ref{after3000}),
$\{\phi_j(x)\}_{j=0}^{\infty}$
is an arbitrary CONS
in $L_2([t, T]),$ and
$\psi_1(\tau),\ldots, \psi_k(\tau)\in L_2([t, T])$.

The equality (\ref{2025may23}) is an analogue
of the formula (5.1) (see \cite{bugh3}, Theorem~5.1)
for the case of a multidimensional Wiener process.

We have

\vspace{-3mm}
$$
T^{k-2r,\hspace{0.2mm}p}_{g_1,g_2,\ldots,g_{2r-1}, g_{2r}}
K(t_{q_1},\ldots,t_{q_{k-2r}})
=
$$

$$
=
\sum\limits_{j_{q_1},\ldots,j_{q_{k-2r}}=0}^p
\Biggl(
\sum\limits_{j_{g_1}, j_{g_3},\ldots ,j_{g_{2r-1}}=0}^p
C_{j_k\ldots j_1}\biggl|_{j_{g_1}=j_{g_2},\ldots, j_{g_{2r-1}}=j_{g_{2r}}}-\Biggr.
$$
$$
\Biggl.-\frac{1}{2^r} \prod\limits_{l=1}^r {\bf 1}_{\{g_{2l}=g_{2l-1}+1\}}
C_{j_k \ldots j_1}\biggl|_{(j_{g_2} j_{g_1})\curvearrowright (\cdot)
\ldots (j_{g_{2r}} j_{g_{2r-1}})\curvearrowright (\cdot),
j_{g_{{}_{1}}}=~j_{g_{{}_{2}}},\ldots, j_{g_{{}_{2r-1}}}=~j_{g_{{}_{2r}}}
}\biggr.\Biggr)\times
$$

\vspace{4mm}

$$
\times\phi_{j_{q_1}}(t_{q_1})\ldots \phi_{j_{q_{k-2r}}}(t_{q_{k-2r}})+
$$

\vspace{-5mm}
$$
+\frac{1}{2^r} \prod\limits_{l=1}^r {\bf 1}_{\{g_{2l}=g_{2l-1}+1\}}
\hspace{-1.5mm}\sum\limits_{j_{q_1},\ldots,j_{q_{k-2r}}=0}^p\hspace{-1mm}
C_{j_k \ldots j_1}\biggl|_{(j_{g_2} j_{g_1})\curvearrowright (\cdot)
\ldots (j_{g_{2r}} j_{g_{2r-1}})\curvearrowright (\cdot),
j_{g_{{}_{1}}}=~j_{g_{{}_{2}}},\ldots, j_{g_{{}_{2r-1}}}=~j_{g_{{}_{2r}}}
}\biggr.
\hspace{-2mm}\times
$$
$$
\times\phi_{j_{q_1}}(t_{q_1})\ldots \phi_{j_{q_{k-2r}}}(t_{q_{k-2r}})\stackrel{\sf def}{=}
$$

\vspace{-2mm}
\begin{equation}
\label{2025may28}
\stackrel{\sf def}{=}F^{(p)}_{g_1,g_2,\ldots,g_{2r-1}, g_{2r}}(t_{q_1},\ldots,t_{q_{k-2r}})+
G^{(p)}_{g_1,g_2,\ldots,g_{2r-1}, g_{2r}}(t_{q_1},\ldots,t_{q_{k-2r}}).
\end{equation}

\vspace{3mm}

Denote
$$
C_{j_k \ldots j_1}\biggl|_{(j_{g_2} j_{g_1})\curvearrowright (\cdot)
\ldots (j_{g_{2r}} j_{g_{2r-1}})\curvearrowright (\cdot),
j_{g_{{}_{1}}}=~j_{g_{{}_{2}}},\ldots, j_{g_{{}_{2r-1}}}=~j_{g_{{}_{2r}}}
}\biggr.\stackrel{\sf def}{=}C_{j_{q_{k-2r}}\ldots j_{q_1}}^{g_1 g_2\ldots g_{2r-1} g_{2r}}.
$$ 

\vspace{2mm}

Using Fubini's Theorem, we obtain (see the derivation of (135) in \cite{2024xxx})
$$
\int\limits_t^T h_{k}(t_k)\ldots \int\limits_t^{t_{l+2}} h_{l+1}(t_{l+1})
\int\limits_t^{t_{l+1}} h_{l}(t_{l})
\int\limits_t^{t_{l}} h_{l-1}(t_{l-1})\ldots
\int\limits_t^{t_2} h_{1}(t_1)
dt_1\ldots 
$$

\vspace{-3mm}
$$
\ldots
dt_{l-1}dt_{l}dt_{l+1}\ldots dt_k=
$$

\vspace{-1mm}
$$
=\int\limits_t^T h_{k}(t_k)\ldots \int\limits_t^{t_{l+2}} h_{l+1}(t_{l+1})
\left(\int\limits_t^{t_{l+1}} h_{l}(t_{l})dt_l\right)
\int\limits_t^{t_{l+1}} h_{l-1}(t_{l-1})\ldots
$$
$$
\ldots 
\int\limits_t^{t_2} h_{1}(t_1)
dt_1\ldots dt_{l-1}dt_{l+1}\ldots dt_k-
$$
$$
-\int\limits_t^T h_{k}(t_k)\ldots \int\limits_t^{t_{l+2}} h_{l+1}(t_{l+1})
\int\limits_t^{t_{l+1}} h_{l-1}(t_{l-1})\left(\int\limits_t^{t_{l-1}} h_{l}(t_{l})dt_l\right)
\int\limits_t^{t_{l-1}} h_{l-2}(t_{l-2})
\ldots
$$
\begin{equation}
\label{july100003}
\ldots
\int\limits_t^{t_2} h_{1}(t_1)
dt_1\ldots dt_{l-2}dt_{l-1}dt_{l+1}\ldots dt_k,
\end{equation}

\vspace{1mm}
\noindent
where $2<l<k-1$ and $h_1(\tau),\ldots,h_k(\tau)\in L_2([t, T]).$ 

By analogy with (\ref{july100003}) we have for $l=k$ (see the derivation of (136) in \cite{2024xxx})
$$
\int\limits_t^{T} h_{l}(t_{l})
\int\limits_t^{t_{l}} h_{l-1}(t_{l-1})\ldots
\int\limits_t^{t_2} h_{1}(t_1)
dt_1\ldots 
dt_{l-1}dt_{l}=
$$
$$
=\left(\int\limits_{t}^{T} h_{l}(t_{l})
dt_l\right)\int\limits_{t}^{T} h_{l-1}(t_{l-1})
\ldots
\int\limits_{t}^{t_2} h_{1}(t_{1})
dt_{1}\ldots dt_{l-1}-
$$
\begin{equation}
\label{july100004}
-\int\limits_{t}^{T}
h_{l-1}(t_{l-1})\left(
\int\limits_{t}^{t_{l-1}} h_{l}(t_{l})dt_l\right)
\int\limits_{t}^{t_{l-1}} h_{l-2}(t_{l-2})\ldots
\int\limits_{t}^{t_2} 
h_{1}(t_{1})
dt_{1}\ldots dt_{l-1}.
\end{equation}

\vspace{2mm}

We will assume that for $l=1$ the transformation 
(\ref{july100003}) is not carried out since
$$
\int\limits_t^{t_2} h_{1}(t_1)
dt_1
$$

\noindent
is the innermost integral on the left-hand side of (\ref{july100003}).

Applying transformations 
(\ref{july100003}), (\ref{july100004}) 
iteratively 
to $C_{j_{q_{k-2r}}\ldots j_{q_1}}^{g_1 g_2\ldots g_{2r-1} g_{2r}}$
for integrations not involving the basis functions
$\phi_{j_{q_1}},\ldots, \phi_{j_{q_{k-2r}}},$ 
we obtain
\begin{equation}
\label{2025may24}
C_{j_{q_{k-2r}}\ldots j_{q_1}}^{g_1 g_2\ldots g_{2r-1}  g_{2r}}=
\sum\limits_{d=1}^{2^{r}}(-1)^{d-1}
\left(\hat C_{j_{q_{k-2r}}\ldots j_{q_1}}^{(d)g_1 g_2\ldots g_{2r-1} g_{2r}}
-
\bar C_{j_{q_{k-2r}}\ldots j_{q_1}}^{(d)g_1 g_2\ldots g_{2r-1} g_{2r}}
\right),
\end{equation}

\vspace{3mm}
\noindent
where some terms in the sum
$$
\sum\limits_{d=1}^{2^{r}}
$$

\noindent
can be identically equal to zero due to
the remark to (\ref{july100003}), (\ref{july100004}).

Using (\ref{2025may24}), we get
$$
G^{(p)}_{g_1,g_2,\ldots,g_{2r-1}, g_{2r}}(t_{q_1},\ldots,t_{q_{k-2r}})=
\frac{1}{2^r} \prod\limits_{l=1}^r {\bf 1}_{\{g_{2l}=g_{2l-1}+1\}}\times
$$
$$
\times
\sum\limits_{d=1}^{2^{r}}(-1)^{d-1}\Biggl(
\sum\limits_{j_{q_1},\ldots,j_{q_{k-2r}}=0}^p
\hat C_{j_{q_{k-2r}}\ldots j_{q_1}}^{(d)g_1 g_2\ldots g_{2r-1} g_{2r}}\cdot
\phi_{j_{q_1}}(t_{q_1})\ldots \phi_{j_{q_{k-2r}}}(t_{q_{k-2r}})-\Biggr.
$$
$$
\Biggl.-
\sum\limits_{j_{q_1},\ldots,j_{q_{k-2r}}=0}^p
\bar C_{j_{q_{k-2r}}\ldots j_{q_1}}^{(d)g_1 g_2\ldots g_{2r-1} g_{2r}}\cdot
\phi_{j_{q_1}}(t_{q_1})\ldots \phi_{j_{q_{k-2r}}}(t_{q_{k-2r}})\Biggr)\ \rightarrow
$$
$$
\rightarrow
\frac{1}{2^r} \prod\limits_{l=1}^r {\bf 1}_{\{g_{2l}=g_{2l-1}+1\}}\times
$$
$$
\times
\sum\limits_{d=1}^{2^{r}}(-1)^{d-1}\biggl(
\hat F^{(d)}_{g_1,g_2,\ldots,g_{2r-1}, g_{2r}}(t_{q_1},\ldots,t_{q_{k-2r}})-
\bar F^{(d)}_{g_1,g_2,\ldots,g_{2r-1}, g_{2r}}(t_{q_1},\ldots,t_{q_{k-2r}})\biggr)\stackrel{\sf def}{=}
$$

\vspace{-2mm}
\begin{equation}
\label{2025may25}
\stackrel{\sf def}{=}
G_{g_1,g_2,\ldots,g_{2r-1}, g_{2r}}(t_{q_1},\ldots,t_{q_{k-2r}})\ \ \hbox{if}\ \ p\to\infty\ \ 
(\hbox{in}\ L_2([t, T]^{k-2r})), 
\end{equation}

\vspace{3mm}
\noindent
where $d=1,\ldots,2^{r},$ $G_{g_1,g_2,\ldots,g_{2r-1}, g_{2r}}(t_{q_1},\ldots,t_{q_{k-2r}}),$ 
$\hat F^{(d)}_{g_1,g_2,\ldots,g_{2r-1}, g_{2r}}(t_{q_1},\ldots,t_{q_{k-2r}}),$
$\bar F^{(d)}_{g_1,g_2,\ldots,g_{2r-1}, g_{2r}}(t_{q_1},\ldots,t_{q_{k-2r}})\in L_2([t, T]^{k-2r}).$

Futhermore, 

\vspace{-4mm}
$$
\bigl\Vert F^{(p)}_{g_1,g_2,\ldots,g_{2r-1}, g_{2r}} \bigr\Vert_{L_2([t, T]^{k-2r})}^2=
$$

\vspace{-2mm}
$$
=\int\limits_{[t, T]^{k-2r}}
\Biggl(
\sum\limits_{j_{q_1},\ldots,j_{q_{k-2r}}=0}^p
\Biggl(
\sum\limits_{j_{g_1}, j_{g_3},\ldots ,j_{g_{2r-1}}=0}^p
C_{j_k\ldots j_1}\biggl|_{j_{g_1}=j_{g_2},\ldots, j_{g_{2r-1}}=j_{g_{2r}}}-\Biggr.\Biggr.
$$

\vspace{-1mm}
$$
\Biggl.-\frac{1}{2^r} \prod\limits_{l=1}^r {\bf 1}_{\{g_{2l}=g_{2l-1}+1\}}
C_{j_k \ldots j_1}\biggl|_{(j_{g_2} j_{g_1})\curvearrowright (\cdot)
\ldots (j_{g_{2r}} j_{g_{2r-1}})\curvearrowright (\cdot),
j_{g_{{}_{1}}}=~j_{g_{{}_{2}}},\ldots, j_{g_{{}_{2r-1}}}=~j_{g_{{}_{2r}}}
}\biggr.\Biggr)\times
$$

\vspace{2mm}
$$
\Biggl.\times\phi_{j_{q_1}}(t_{q_1})\ldots \phi_{j_{q_{k-2r}}}(t_{q_{k-2r}})\Biggr)^2
dt_{q_1}\ldots dt_{q_{k-2r}}=
$$

\vspace{3mm}
$$
=
\sum\limits_{j_{q_1},\ldots,j_{q_{k-2r}}=0}^p
\Biggl(
\sum\limits_{j_{g_1}, j_{g_3},\ldots ,j_{g_{2r-1}}=0}^p
C_{j_k\ldots j_1}\biggl|_{j_{g_1}=j_{g_2},\ldots, j_{g_{2r-1}}=j_{g_{2r}}}-\Biggr.\Biggr.
$$

\vspace{-1mm}
$$
\Biggl.-\frac{1}{2^r} \prod\limits_{l=1}^r {\bf 1}_{\{g_{2l}=g_{2l-1}+1\}}
C_{j_k \ldots j_1}\biggl|_{(j_{g_2} j_{g_1})\curvearrowright (\cdot)
\ldots (j_{g_{2r}} j_{g_{2r-1}})\curvearrowright (\cdot),
j_{g_{{}_{1}}}=~j_{g_{{}_{2}}},\ldots, j_{g_{{}_{2r-1}}}=~j_{g_{{}_{2r}}}
}\biggr.\Biggr)^2.
$$

\vspace{5mm}

This means that the condition (\ref{july700000}) is equivalent to 
\begin{equation}
\label{2025may26}
\lim\limits_{p\to\infty}
\bigl\Vert F^{(p)}_{g_1,g_2,\ldots,g_{2r-1}, g_{2r}} \bigr\Vert_{L_2([t, T]^{k-2r})}^2=0.
\end{equation}

\vspace{1mm}

Suppose that the condition (\ref{july700000}) (or (\ref{2025may26}))
is fulfilled. 
Applying (\ref{2025may28}), (\ref{2025may25}) and (\ref{2025may26}), we obtain
$$
\bigl\Vert T^{k-2r,\hspace{0.2mm}p}_{g_1,g_2,\ldots,g_{2r-1}, g_{2r}}
K- G_{g_1,g_2,\ldots,g_{2r-1}, g_{2r}} \bigr\Vert_{L_2([t, T]^{k-2r})}=
$$
$$
=\bigl\Vert F^{(p)}_{g_1,g_2,\ldots,g_{2r-1}, g_{2r}} + 
G^{(p)}_{g_1,g_2,\ldots,g_{2r-1}, g_{2r}} - G_{g_1,g_2,\ldots,g_{2r-1}, g_{2r}} 
\bigr\Vert_{L_2([t, T]^{k-2r})}\le
$$

$$
\le\bigl\Vert F^{(p)}_{g_1,g_2,\ldots,g_{2r-1}, g_{2r}} \bigr\Vert_{L_2([t, T]^{k-2r})} +
$$

$$
+ 
\bigl\Vert G^{(p)}_{g_1,g_2,\ldots,g_{2r-1}, g_{2r}} - 
G_{g_1,g_2,\ldots,g_{2r-1}, g_{2r}} \bigr\Vert_{L_2([t, T]^{k-2r})}\ \to 0
$$

\vspace{2mm}
\noindent
if $p\to\infty.$

Thus, the limiting trace $T^{k-2r}_{g_1,g_2,\ldots,g_{2r-1}, g_{2r}}
K(t_{q_1},\ldots,t_{q_{k-2r}})$
exists under the condition (\ref{july700000})
for all $r=1,2,\ldots,$ $[k/2]$ and 
for all possible $g_1,g_2,\ldots,g_{2r-1},g_{2r}$ 
(see (\ref{leto5007after})), i.e.

\vspace{-4mm}
$$
T^{k-2r}_{g_1,g_2,\ldots,g_{2r-1}, g_{2r}}K(t_{q_1},\ldots,t_{q_{k-2r}})
=
\lim\limits_{p\to\infty}T^{k-2r,\hspace{0.2mm}p}_{g_1,g_2,\ldots,g_{2r-1}, g_{2r}}
K(t_{q_1},\ldots,t_{q_{k-2r}})
=
$$

\vspace{-1mm}
\begin{equation}
\label{october20251}
=G_{g_1,g_2,\ldots,g_{2r-1}, g_{2r}}(t_{q_1},\ldots,t_{q_{k-2r}})
\end{equation}

\vspace{4mm}
\noindent
in $L_2([t, T]^{k-2r}),$ where 
$G_{g_1,g_2,\ldots,g_{2r-1}, g_{2r}}(t_{q_1},\ldots,t_{q_{k-2r}})$ is defined by (\ref{2025may25}).

Further, we assume that the condition (\ref{october20251}) is fulfilled.
Using (\ref{october20251}) and (\ref{2025may28}), (\ref{2025may25}), we obtain 

\vspace{-3mm}
$$
\bigl\Vert F^{(p)}_{g_1,g_2,\ldots,g_{2r-1}, g_{2r}}\bigr\Vert_{L_2([t, T]^{k-2r})}=
$$

$$
=\bigl\Vert T^{k-2r,\hspace{0.2mm}p}_{g_1,g_2,\ldots,g_{2r-1}, g_{2r}}
K- G^{(p)}_{g_1,g_2,\ldots,g_{2r-1}, g_{2r}} \bigr\Vert_{L_2([t, T]^{k-2r})}=
$$

$$
=\bigl\Vert T^{k-2r,\hspace{0.2mm}p}_{g_1,g_2,\ldots,g_{2r-1}, g_{2r}}
K- G_{g_1,g_2,\ldots,g_{2r-1}, g_{2r}}+
G_{g_1,g_2,\ldots,g_{2r-1}, g_{2r}}
-
G^{(p)}_{g_1,g_2,\ldots,g_{2r-1}, g_{2r}} \bigr\Vert_{L_2([t, T]^{k-2r})}
$$

$$
\le\bigl\Vert T^{k-2r,\hspace{0.2mm}p}_{g_1,g_2,\ldots,g_{2r-1}, g_{2r}}
K- G_{g_1,g_2,\ldots,g_{2r-1}, g_{2r}} \bigr\Vert_{L_2([t, T]^{k-2r})} + 
$$

$$
+
\bigl\Vert G_{g_1,g_2,\ldots,g_{2r-1}, g_{2r}}
-
G^{(p)}_{g_1,g_2,\ldots,g_{2r-1}, g_{2r}} \bigr\Vert_{L_2([t, T]^{k-2r})}\ \to 0
$$

\vspace{3mm}
\noindent
if $p\to\infty.$

Thus, the condition (\ref{2025may26}) is satisfied
under the condition (\ref{october20251}).
Note that if $T^{k-2r}_{g_1,g_2,\ldots,g_{2r-1}, g_{2r}}K(t_{q_1},\ldots,t_{q_{k-2r}})$
exists, then it will be equal to 
$G_{g_1,g_2,\ldots,g_{2r-1}, g_{2r}}(t_{q_1},\ldots,t_{q_{k-2r}}),$ which is defined by (\ref{2025may25})
(see Sect.~3.2 for details).

\subsection{New Representations of the Hu--Meyer Formulas for the Case
of a Multidimensional Wiener Process. Connection with Theorems~4 and 7}

Suppose that $\{\phi_j(x)\}_{j=0}^{\infty}$
is an arbitrary CONS of functions
in $L_2([t, T])$ and
$\Phi(t_1,\ldots,t_k)\in L_2([t, T]^k).$

Let us generalize the definition of the limiting trace
from the previous section.

We define the 
$r$th limiting trace of the function $\Phi(t_1,\ldots,t_k)\in L_2([t, T]^k)$
by the following expression

\vspace{-3mm}
$$
T^{k-2r}_{g_1,g_2,\ldots,g_{2r-1}, g_{2r}}\Phi (t_{q_1},\ldots,t_{q_{k-2r}})\stackrel{\sf def}{=}
$$

\vspace{-2mm}
\begin{equation}
\label{novem2026xxx1122}
\stackrel{\sf def}{=}
\lim\limits_{p\to\infty}T^{k-2r,\hspace{0.2mm}p}_{g_1,g_2,\ldots,g_{2r-1}, g_{2r}}
\Phi (t_{q_1},\ldots,t_{q_{k-2r}})
\end{equation}

\vspace{1mm}
\noindent
in $L_2([t, T]^{k-2r})$ (here and further $L_2([t, T]^{0})\stackrel{\sf def}{=}{\bf R}$), where

\vspace{-3mm}
$$
T^{k-2r,\hspace{0.2mm}p}_{g_1,g_2,\ldots,g_{2r-1}, g_{2r}}\Phi (t_{q_1},\ldots,t_{q_{k-2r}})=
$$

\vspace{-4mm}
$$
=
\sum\limits_{j_{q_1},\ldots,j_{q_{k-2r}}=0}^p
\sum\limits_{j_{g_1}, j_{g_3},\ldots ,j_{g_{2r-1}}=0}^p
C_{j_k\ldots j_1}\biggl|_{j_{g_1}=j_{g_2},\ldots, j_{g_{2r-1}}=j_{g_{2r}}}
\times
$$

\vspace{1mm}
\begin{equation}
\label{trace11111111}
\times
\phi_{j_{q_1}}(t_{q_1})\ldots \phi_{j_{q_{k-2r}}}(t_{q_{k-2r}}),
\end{equation}

\vspace{3mm}
\noindent
where
$r=1,2,\ldots,[k/2]$ and 
$\{g_1,g_2,\ldots,g_{2r-1},g_{2r},
q_1,\ldots,q_{k-2r}\}=\{1,2,\ldots,k\}$ 
(see (\ref{leto5007after})),

\vspace{-3mm}
\begin{equation}
\label{2025may30}
C_{j_k\ldots j_1}=\int\limits_{[t,T]^k}
\Phi(t_1,\ldots,t_k)\prod_{l=1}^{k}\phi_{j_l}(t_l)dt_1\ldots dt_k
\end{equation}

\vspace{2mm}
\noindent
is the Fourier coefficient.
Also we will write
$T^{k-2r}_{g_1,g_2,\ldots,g_{2r-1}, g_{2r}}\Phi (t_{q_1},\ldots,t_{q_{k-2r}})=\Phi(t_1,\ldots,t_k)$ 
and 
$T^{k-2r,\hspace{0.2mm}p}_{g_1,g_2,\ldots,g_{2r-1}, g_{2r}}
\Phi (t_{q_1},\ldots,t_{q_{k-2r}})=\Phi_p(t_1,\ldots,t_k)$
for $r=0$, where

\vspace{-2mm}
$$
\Phi_p(t_1,\ldots,t_k)=
\sum\limits_{j_{1},\ldots,j_{k}=0}^p
C_{j_k\ldots j_1}
\phi_{j_{1}}(t_{1})\ldots \phi_{j_k}(t_k).
$$

\vspace{2mm}

Let us consider a prelimit variant of the formula (\ref{2025may23})
for the case $p_1=\ldots=p_k=p$ (see the proof in \cite{2018a}, Sect.~2.18) and replace 
$K(t_1,\ldots,t_k)$ with $\Phi(t_1,\ldots,t_k)$ in it. Thus, we have

\vspace{-1mm}
$$
\sum_{j_1,\ldots,j_k=0}^{p}
C_{j_k\ldots j_1}
\prod_{l=1}^k \zeta_{j_l}^{(i_l)}=
\sum_{j_1,\ldots,j_k=0}^{p}
C_{j_k\ldots j_1}
J'[\phi_{j_1}\ldots \phi_{j_k}]_{T,t}^{(i_1\ldots i_k)}+
$$

\vspace{2mm}
$$
+
\sum\limits_{r=1}^{[k/2]}
\sum_{\stackrel{(\{\{g_1, g_2\}, \ldots, 
\{g_{2r-1}, g_{2r}\}\}, \{q_1, \ldots, q_{k-2r}\})}
{{}_{\{g_1, g_2, \ldots, 
g_{2r-1}, g_{2r}, q_1, \ldots, q_{k-2r}\}=\{1, 2, \ldots, k\}}}}
\prod\limits_{s=1}^r
{\bf 1}_{\{i_{g_{{}_{2s-1}}}=~i_{g_{{}_{2s}}}\ne 0\}}
\sum_{j_1,\ldots,j_k=0}^{p}
C_{j_k\ldots j_1}\times
$$

\vspace{2mm}
\begin{equation}
\label{2025may31}
\times{\bf 1}_{\{j_{g_{{}_{2s-1}}}=~j_{g_{{}_{2s}}}\}}
J'[\phi_{j_{q_1}}\ldots \phi_{j_{q_{k-2r}}}]_{T,t}^{(i_{q_1}\ldots i_{q_{k-2r}})}\ \ \ \hbox{w.~p.~1,}
\end{equation}

\vspace{3mm}
\noindent
where $k\ge 2,$ $J'[\phi_{j_1}\ldots \phi_{j_k}]_{T,t}^{(i_1\ldots i_k)},$
$J'[\phi_{j_{q_1}}\ldots \phi_{j_{q_{k-2r}}}]_{T,t}^{(i_{q_1}\ldots i_{q_{k-2r}})}$
are multiple Wie\-ner sto\-chas\-tic integrals 
defined as in \cite{ito1951} but for the case of a multidimensional Wiener
process (see \cite{2018a}, Sect.~1.11 for details),
$J'[\phi_{j_{q_1}}\ldots \phi_{j_{q_{k-2r}}}]_{T,t}^{(i_{q_1}\ldots i_{q_{k-2r}})}
\stackrel{\sf def}{=}1$ for $k=2r,$ and $C_{j_k\ldots j_1}$
is defined by (\ref{2025may30}).

Note that sometimes the multiple Stratonovich stochastic integral for
$\Phi(t_1,\ldots,t_k)\in L_2([t, T]^k)$
is defined as the limit in probability as $p\to\infty$ of the following expression
(the case of a scalar standard Wiener process)
$$
\sum_{j_1,\ldots,j_k=0}^{p}
C_{j_k\ldots j_1}
\prod_{l=1}^k \zeta_{j_l}^{(i)},
$$

\noindent
where $i_1=\ldots=i_k=i$ (see, for example,  Definition~5.9 in \cite{HuHu}).

By analogy with \cite{HuHu}, we define the multiple Stratonovich
stochastic integral
for $\Phi(t_1,\ldots,t_k)\in L_2([t, T]^k)$
(the case of a multidimensional Wiener process)
as the following mean-square limit

\vspace{-2mm}
\begin{equation}
\label{novem2026xxx4}
\hat J^{S}[\Phi]_{T,t}^{(i_1\ldots i_k)}=
\hbox{\vtop{\offinterlineskip\halign{
\hfil#\hfil\cr
{\rm l.i.m.}\cr
$\stackrel{}{{}_{p\to \infty}}$\cr
}} }
\hat J^{S}_p[\Phi]_{T,t}^{(i_1\ldots i_k)},
\end{equation}

\noindent
where
\begin{equation}
\label{novem2026xxx5}
\hat J^{S}_p[\Phi]_{T,t}^{(i_1\ldots i_k)}=
\sum_{j_1,\ldots,j_k=0}^{p}
C_{j_k\ldots j_1}
\prod_{l=1}^k \zeta_{j_l}^{(i_l)},
\end{equation}

\vspace{2mm}
\noindent
$C_{j_k\ldots j_1}$ is the Fourier coefficient
defined by (\ref{2025may30}).

Let us rewrite (\ref{2025may31}) in the form

\vspace{-2mm}
$$
\hat J^{S}_p[\Phi]_{T,t}^{(i_1\ldots i_k)}
=
J'[\Phi_p]_{T,t}^{(i_1\ldots i_k)}+
$$

\vspace{-2mm}
$$
+
\sum\limits_{r=1}^{[k/2]}
\sum_{\stackrel{(\{\{g_1, g_2\}, \ldots, 
\{g_{2r-1}, g_{2r}\}\}, \{q_1, \ldots, q_{k-2r}\})}
{{}_{\{g_1, g_2, \ldots, 
g_{2r-1}, g_{2r}, q_1, \ldots, q_{k-2r}\}=\{1, 2, \ldots, k\}}}}
\prod\limits_{s=1}^r
{\bf 1}_{\{i_{g_{{}_{2s-1}}}=~i_{g_{{}_{2s}}}\ne 0\}}\times
$$

\vspace{5mm}
\begin{equation}
\label{2026may1}
\times
J'\hspace{-1mm}\left[T^{k-2r,\hspace{0.2mm}p}_{g_1,g_2,\ldots,g_{2r-1}, 
g_{2r}}\Phi\right]_{T,t}^{(i_{q_1}\ldots i_{q_{k-2r}})}\ \ \ \hbox{w.~p.~1,}
\end{equation}

\vspace{5mm}
\noindent
where $J'\hspace{-1mm}\left[T^{k-2r,\hspace{0.2mm}p}_{g_1,g_2,\ldots,g_{2r-1}, 
g_{2r}}\Phi\right]_{T,t}^{(i_{q_1}\ldots i_{q_{k-2r}})}\stackrel{\sf def}{=}
T^{k-2r,\hspace{0.2mm}p}_{g_1,g_2,\ldots,g_{2r-1}, 
g_{2r}}\Phi$ for $k=2r.$

Assume that all limiting traces
$T^{k-2r}_{g_1,g_2,\ldots,g_{2r-1}, g_{2r}}\Phi (t_{q_1},\ldots,t_{q_{k-2r}})$
(for all $r=1,2,\ldots,$ $[k/2]$ and for all possible
$g_1,g_2,\ldots,g_{2r-1},g_{2r}$ (see {\rm (\ref{leto5007after})))
exist.

Recall the well-known property
of the multiple Wiener stochastic integral 

\vspace{-3mm}
\begin{equation}
\label{2026may2}
{\sf M}\left\{\left(J'[\Phi]_{T,t}^{(i_1\ldots i_k)}\right)^2\right\}\le
C_k
\int\limits_{[t,T]^k}
\Phi^2(t_1,\ldots,t_k)dt_1\ldots dt_k,
\end{equation}

\vspace{2mm}
\noindent
where $\Phi(t_1,\ldots,t_k)\in L_2([t, T]^k)$ and $C_k$ is a constant.

Using (\ref{2026may1}), (\ref{2026may2}) and the existence of limiting traces, we have

\vspace{-3mm}
$$
{\sf M}\Biggl\{\Biggl(\hat J^{S}_p[\Phi]_{T,t}^{(i_1\ldots i_k)}
-
J'[\Phi]_{T,t}^{(i_1\ldots i_k)}-\Biggr.\Biggr.
$$

\vspace{-2mm}

$$
-
\sum\limits_{r=1}^{[k/2]}
\sum_{\stackrel{(\{\{g_1, g_2\}, \ldots, 
\{g_{2r-1}, g_{2r}\}\}, \{q_1, \ldots, q_{k-2r}\})}
{{}_{\{g_1, g_2, \ldots, 
g_{2r-1}, g_{2r}, q_1, \ldots, q_{k-2r}\}=\{1, 2, \ldots, k\}}}}
\prod\limits_{s=1}^r
{\bf 1}_{\{i_{g_{{}_{2s-1}}}=~i_{g_{{}_{2s}}}\ne 0\}}\times
$$

$$
\Biggl.\Biggl.\times
J'\hspace{-1mm}\left[T^{k-2r}_{g_1,g_2,\ldots,g_{2r-1}, 
g_{2r}}\Phi\right]_{T,t}^{(i_{q_1}\ldots i_{q_{k-2r}})}\Biggr)^2\Biggr\}=
$$

\vspace{-1mm}
$$
={\sf M}\Biggl\{\Biggl(
J'[\Phi_p-\Phi]_{T,t}^{(i_1\ldots i_k)}+\Biggr.\Biggr.
$$
$$
+
\sum\limits_{r=1}^{[k/2]}
\sum_{\stackrel{(\{\{g_1, g_2\}, \ldots, 
\{g_{2r-1}, g_{2r}\}\}, \{q_1, \ldots, q_{k-2r}\})}
{{}_{\{g_1, g_2, \ldots, 
g_{2r-1}, g_{2r}, q_1, \ldots, q_{k-2r}\}=\{1, 2, \ldots, k\}}}}
\prod\limits_{s=1}^r
{\bf 1}_{\{i_{g_{{}_{2s-1}}}=~i_{g_{{}_{2s}}}\ne 0\}}\times
$$

$$
\Biggl.\Biggl.\times
J'\hspace{-1mm}\left[T^{k-2r,\hspace{0.2mm}p}_{g_1,g_2,\ldots,g_{2r-1}, 
g_{2r}}\Phi-
T^{k-2r}_{g_1,g_2,\ldots,g_{2r-1}, 
g_{2r}}\Phi
\right]_{T,t}^{(i_{q_1}\ldots i_{q_{k-2r}})}\Biggr)^2\Biggr\}\le
$$

\vspace{1mm}

$$
\le C_k'\Biggl(
\bigl\Vert\Phi_p-\Phi\bigr\Vert_{L_2([t, T]^k)}^2+\Biggr.
$$

\vspace{-3mm}
$$
+
\sum\limits_{r=1}^{[k/2]}
\sum_{\stackrel{(\{\{g_1, g_2\}, \ldots, 
\{g_{2r-1}, g_{2r}\}\}, \{q_1, \ldots, q_{k-2r}\})}
{{}_{\{g_1, g_2, \ldots, 
g_{2r-1}, g_{2r}, q_1, \ldots, q_{k-2r}\}=\{1, 2, \ldots, k\}}}}
\prod\limits_{s=1}^r
{\bf 1}_{\{i_{g_{{}_{2s-1}}}=~i_{g_{{}_{2s}}}\ne 0\}}\times
$$

\begin{equation}
\label{2026may5}
\Biggl.\times
\left\Vert T^{k-2r,\hspace{0.2mm}p}_{g_1,g_2,\ldots,g_{2r-1}, 
g_{2r}}\Phi-
T^{k-2r}_{g_1,g_2,\ldots,g_{2r-1}, 
g_{2r}}\Phi\right\Vert_{L_2([t, T]^{k-2r})}^2\Biggr)\ \to\ 0
\end{equation}

\vspace{2mm}
\noindent
if $p\to\infty,$ 
where $J'\hspace{-1mm}\left[T^{k-2r}_{g_1,g_2,\ldots,g_{2r-1}, 
g_{2r}}\Phi\right]_{T,t}^{(i_{q_1}\ldots i_{q_{k-2r}})}\stackrel{\sf def}{=}
T^{k-2r}_{g_1,g_2,\ldots,g_{2r-1}, 
g_{2r}}\Phi$ for $k=2r,$
and $C_k'$ is a constant.

Applying (\ref{2026may5}), we obtain the following 
new representation of the Hu--Meyer formula for the case
of a multidimensional Wiener process

\vspace{-2mm}
$$
\hat J^{S}[\Phi]_{T,t}^{(i_1\ldots i_k)}
=
J'[\Phi]_{T,t}^{(i_1\ldots i_k)}+
$$

\vspace{-3mm}
$$
+
\sum\limits_{r=1}^{[k/2]}
\sum_{\stackrel{(\{\{g_1, g_2\}, \ldots, 
\{g_{2r-1}, g_{2r}\}\}, \{q_1, \ldots, q_{k-2r}\})}
{{}_{\{g_1, g_2, \ldots, 
g_{2r-1}, g_{2r}, q_1, \ldots, q_{k-2r}\}=\{1, 2, \ldots, k\}}}}
\prod\limits_{s=1}^r
{\bf 1}_{\{i_{g_{{}_{2s-1}}}=~i_{g_{{}_{2s}}}\ne 0\}}\times
$$

\vspace{4mm}
\begin{equation}
\label{2026may10}
\times
J'\hspace{-1mm}\left[T^{k-2r}_{g_1,g_2,\ldots,g_{2r-1}, 
g_{2r}}\Phi\right]_{T,t}^{(i_{q_1}\ldots i_{q_{k-2r}})}\ \ \ \hbox{w.~p.~1.}
\end{equation}

\vspace{4mm}

The equality (\ref{2026may10}) is consistent with Theorem~5.1 \cite{bugh3}
(the case of a scalar Wiener process).

Further, let us obtain the inverse version of the Hu--Meyer
formula (\ref{2026may10}), i.e. a formula expressing
the multiple Wiener stochastic integral through the sum
of multiple Stratonovich stochastic integrals.

Let us give the following definition of the limiting trace.
We define the 
$r$th limiting trace of the function $\Phi(t_1,\ldots,t_k)\in L_2([t, T]^k)$
by the following expression
\begin{equation}
\label{soglas123sd}
\tilde T^{k-2r}_{g_1,g_2,\ldots,g_{2r-1}, g_{2r}}
\Phi (t_{q_1},\ldots,t_{q_{k-2r}})\stackrel{\sf def}{=}
\lim\limits_{p\to\infty}\tilde T^{k-2r,\hspace{0.2mm}p}_{g_1,g_2,\ldots,g_{2r-1}, g_{2r}}
\Phi (t_{q_1},\ldots,t_{q_{k-2r}})
\end{equation}

\noindent
in $L_2([t, T]^{k-2r})$, where

\vspace{-3mm}
$$
\tilde T^{k-2r,\hspace{0.2mm}p}_{g_1,g_2,\ldots,g_{2r-1}, g_{2r}}\Phi (t_{q_1},\ldots,t_{q_{k-2r}})=
$$

\vspace{-6mm}
$$
=\sum\limits_{j_{q_1},\ldots,j_{q_{k-2r}}=0}^p
\sum\limits_{j_{g_1}, j_{g_3},\ldots ,j_{g_{2r-1}}=0}^{\infty}
C_{j_k\ldots j_1}\biggl|_{j_{g_1}=j_{g_2},\ldots, j_{g_{2r-1}}=j_{g_{2r}}}
\cdot\hspace{0.7mm} \phi_{j_{q_1}}(t_{q_1})\ldots \phi_{j_{q_{k-2r}}}(t_{q_{k-2r}})=
$$
         
\vspace{1mm}
\begin{equation}
\label{october202620}
=
\sum\limits_{j_{q_1},\ldots,j_{q_{k-2r}}=0}^p
\tilde C_{j_{q_{k-2r}}\ldots j_{q_1}}
\cdot\hspace{0.7mm}\phi_{j_{q_1}}(t_{q_1})\ldots \phi_{j_{q_{k-2r}}}(t_{q_{k-2r}}),
\end{equation}

\vspace{2mm}
\noindent
where
$$
\tilde C_{j_{q_{k-2r}}\ldots j_{q_1}}=
\sum\limits_{j_{g_1}, j_{g_3},\ldots ,j_{g_{2r-1}}=0}^{\infty}
C_{j_k\ldots j_1}\biggl|_{j_{g_1}=j_{g_2},\ldots, j_{g_{2r-1}}=j_{g_{2r}}}=
$$

\vspace{-1mm}
\begin{equation}
\label{october202621}
=\lim\limits_{p\to\infty}\sum\limits_{j_{g_1}, j_{g_3},\ldots ,j_{g_{2r-1}}=0}^{p}
C_{j_k\ldots j_1}\biggl|_{j_{g_1}=j_{g_2},\ldots, j_{g_{2r-1}}=j_{g_{2r}}},
\end{equation}

\vspace{3mm}
\noindent
where $r=1,2,\ldots,[k/2],$
$\{g_1,g_2,\ldots,g_{2r-1},g_{2r},
q_1,\ldots,q_{k-2r}\}=\{1,2,\ldots,k\}$ 
(see (\ref{leto5007after})),
$C_{j_k\ldots j_1}$ has the form
(\ref{2025may30}), 
$\tilde T^{k-2r}_{g_1,g_2,\ldots,g_{2r-1}, g_{2r}}\Phi (t_{q_1},\ldots,t_{q_{k-2r}})=\Phi(t_1,\ldots,t_k)$ 
and 
$\tilde T^{k-2r,\hspace{0.2mm}p}_{g_1,g_2,\ldots,g_{2r-1}, g_{2r}}\Phi (t_{q_1},\ldots,t_{q_{k-2r}})
=\Phi_p(t_1,\ldots,t_k)$
for $r=0$, where

\vspace{-3mm}
$$
\Phi_p(t_1,\ldots,t_k)=
\sum\limits_{j_{1},\ldots,j_{k}=0}^p
C_{j_k\ldots j_1}
\phi_{j_{1}}(t_{1})\ldots \phi_{j_k}(t_k).
$$

\vspace{2mm}

Suppose that all limiting traces
$T^{k-2r}_{g_1,g_2,\ldots,g_{2r-1}, g_{2r}}\Phi (t_{q_1},\ldots,t_{q_{k-2r}})$
exist
for all $r=1,2,\ldots,$ $[k/2]$ and for all possible
$g_1,g_2,$ $\ldots,g_{2r-1},g_{2r}$ (see {\rm (\ref{leto5007after})).

Let us fix
$j'_{q_1}, \ldots, j'_{q_{k-2r}}$
Applying 
the inequality 
of Cauchy--Bu\-ny\-a\-kov\-sky, we obtain
$$
\lim\limits_{p\to\infty}
\biggl\langle T^{k-2r,\hspace{0.2mm}p}_{g_1,g_2,\ldots,g_{2r-1}, g_{2r}}
\Phi (t_{q_1},\ldots,t_{q_{k-2r}})\ ,\  
\phi_{j'_{q_1}}(t_{q_1})\ldots \phi_{j'_{q_{k-2r}}}(t_{q_{k-2r}}) \biggr\rangle_{L_2([t, T]^{k-2r})}=
$$
$$
=
\biggl\langle T^{k-2r}_{g_1,g_2,\ldots,g_{2r-1}, g_{2r}}\Phi (t_{q_1},\ldots,t_{q_{k-2r}})\ ,\  
\phi_{j'_{q_1}}(t_{q_1})\ldots \phi_{j'_{q_{k-2r}}}(t_{q_{k-2r}}) \biggr\rangle_{L_2([t, T]^{k-2r})},
$$

\vspace{2mm}
\noindent
where 
$T^{k-2r,\hspace{0.2mm}p}_{g_1,g_2,\ldots,g_{2r-1}, g_{2r}}\Phi (t_{q_1},\ldots,t_{q_{k-2r}})$
is defined by (\ref{trace11111111}).

Let $p\ge \max\left\{j'_{q_1}, \ldots, j'_{q_{k-2r}}\right\}.$ Then
$$
\biggl\langle T^{k-2r,\hspace{0.2mm}p}_{g_1,g_2,\ldots,g_{2r-1}, g_{2r}}
\Phi (t_{q_1},\ldots,t_{q_{k-2r}})\ ,\  
\phi_{j'_{q_1}}(t_{q_1})\ldots \phi_{j'_{q_{k-2r}}}(t_{q_{k-2r}}) \biggr\rangle_{L_2([t, T]^{k-2r})}=
$$
$$
=\sum\limits_{j_{q_1},\ldots,j_{q_{k-2r}}=0}^p
\sum\limits_{j_{g_1}, j_{g_3},\ldots ,j_{g_{2r-1}}=0}^p
C_{j_k\ldots j_1}\biggl|_{j_{g_1}=j_{g_2},\ldots, j_{g_{2r-1}}=j_{g_{2r}}}\times
$$

\vspace{-1mm}
$$
\times
\biggl\langle\phi_{j_{q_1}}(t_{q_1})\ldots \phi_{j_{q_{k-2r}}}(t_{q_{k-2r}})\ ,\ 
\phi_{j'_{q_1}}(t_{q_1})\ldots \phi_{j'_{q_{k-2r}}}(t_{q_{k-2r}}) \biggr\rangle_{L_2([t, T]^{k-2r})}=
$$

\vspace{-2mm}
$$
=
\sum\limits_{j_{g_1}, j_{g_3},\ldots ,j_{g_{2r-1}}=0}^{p}
C_{j_k\ldots j_1}\biggl|_{j_{g_1}=j_{g_2},\ldots, j_{g_{2r-1}}=j_{g_{2r}},
j_{q_1}=j'_{q_1},\ldots, j_{q_{k-2r}}=j'_{q_{k-2r}}}.
$$

\vspace{4mm}

This means that
$$
\lim\limits_{p\to\infty}\sum\limits_{j_{g_1}, j_{g_3},\ldots ,j_{g_{2r-1}}=0}^{p}
C_{j_k\ldots j_1}\biggl|_{j_{g_1}=j_{g_2},\ldots, j_{g_{2r-1}}=j_{g_{2r}}}=
$$
$$
=\sum\limits_{j_{g_1}, j_{g_3},\ldots ,j_{g_{2r-1}}=0}^{\infty}
C_{j_k\ldots j_1}\biggl|_{j_{g_1}=j_{g_2},\ldots, j_{g_{2r-1}}=j_{g_{2r}}}
$$

\vspace{3mm}
\noindent
is a Fourier coefficient of
$T^{k-2r}_{g_1,g_2,\ldots,g_{2r-1}, g_{2r}}\Phi (t_{q_1},\ldots,t_{q_{k-2r}}),$
i.e. 
\begin{equation}
\label{req1234}
T^{k-2r}_{g_1,g_2,\ldots,g_{2r-1}, g_{2r}}\Phi (t_{q_1},\ldots,t_{q_{k-2r}})=
\tilde T^{k-2r}_{g_1,g_2,\ldots,g_{2r-1}, g_{2r}}\Phi (t_{q_1},\ldots,t_{q_{k-2r}})
\end{equation}

\vspace{3.5mm}
\noindent
in $L_2([t, T]^{k-2r})$.
Then, using the equality (\ref{req1234})
and (\ref{2026may10}), we obtain

\vspace{-2mm}
$$
J'[\Phi]_{T,t}^{(i_1\ldots i_k)}=\hat J^{S}[\Phi]_{T,t}^{(i_1\ldots i_k)}-
$$

\vspace{-1mm}
$$
-
\sum\limits_{r=1}^{[k/2]}
\sum_{\stackrel{(\{\{g_1, g_2\}, \ldots, 
\{g_{2r-1}, g_{2r}\}\}, \{q_1, \ldots, q_{k-2r}\})}
{{}_{\{g_1, g_2, \ldots, 
g_{2r-1}, g_{2r}, q_1, \ldots, q_{k-2r}\}=\{1, 2, \ldots, k\}}}}
\prod\limits_{s=1}^r
{\bf 1}_{\{i_{g_{{}_{2s-1}}}=~i_{g_{{}_{2s}}}\ne 0\}}\times
$$

\vspace{5mm}
\begin{equation}
\label{2026may10sss}
\times
J'\hspace{-1mm}\left[\tilde T^{k-2r}_{g_1,g_2,\ldots,g_{2r-1}, 
g_{2r}}\Phi\right]_{T,t}^{(i_{q_1}\ldots i_{q_{k-2r}})}\ \ \ \hbox{w.~p.~1.}
\end{equation}

\vspace{5mm}

Since $\tilde T^{k-2r}_{g_1,g_2,\ldots,g_{2r-1}, 
g_{2r}}\Phi(t_{q_1},\ldots,t_{q_{k-2r}})\in L_2([t,T]^{k-2r})$,
then the multiple Stra\-to\-no\-vich stochastic integral 
$\hat J^{S}[\tilde T^{k-2r}_{g_1,g_2,\ldots,g_{2r-1}, g_{2r}}\Phi]_{T,t}^{(i_{q_1}\ldots i_{q_{k-2r}})}$
is defined according to (\ref{novem2026xxx4}), (\ref{novem2026xxx5}).
If we replace $\Phi$ with $\tilde T^{k-2r}_{g_1,g_2,\ldots,g_{2r-1}, g_{2r}}\Phi$
in (\ref{2026may10sss}), then the integral 
$\hat J^{S}[\tilde T^{k-2r}_{g_1,g_2,\ldots,g_{2r-1}, g_{2r}}\Phi]_{T,t}^{(i_{q_1}\ldots i_{q_{k-2r}})}$
will exist according the obtained version of the formula (\ref{2026may10sss})
under some additional conditions (see below).

When applying the formula (\ref{2026may10sss})
iteratively, we will obviously encounter with iterative 
application of traces $\tilde T^{k-2l}_{g_1,g_2,\ldots,g_{2l-1}, 
g_{2l}}\Phi(t_{q_1},\ldots,t_{q_{k-2l}})$ of different orders
$l.$ Let us assume that
$$
\tilde T^{k-2m}_{k_1,k_2,\ldots,k_{2m-1}, k_{2m}}
\tilde T^{k-2l}_{g_1,g_2,\ldots,g_{2l-1},g_{2l}}
\Phi(t_{q_1},\ldots,t_{q_{k-2(l+m)}})=
$$
\begin{equation}
\label{sogl111ds}
=
\tilde T^{k-2(l+m)}_{k_1,k_2,\ldots,k_{2m-1}, k_{2m},g_1,g_2,\ldots,g_{2l-1},g_{2l}}
\Phi(t_{q_1},\ldots,t_{q_{k-2(l+m)}})
\end{equation}

\vspace{1mm}
\noindent
in $L_2([t,T]^{k-2(l+m)}),$ where $l=0,1,2,\ldots,[k/2],$ $m=0,1,2,\ldots,[(k-2l)/2].$

We also suppose that conditions similar to (\ref{sogl111ds})
will be satisfied for the superposition of 3 or more traces
defined by the equality (\ref{soglas123sd}).

Note that the listed conditions will be satisfied if
the following condition is fulfilled.

{\it We will say that
condition {\rm $(\star)$} is satisfied if
for all $r=1,2,\ldots,[k/2]$ and for all possible $g_1,g_2,\ldots,$ $g_{2r-1},g_{2r}$ 
the sum of the series

\newpage
\noindent
$$
\sum\limits_{j_{g_1}, j_{g_3},\ldots ,j_{g_{2r-1}}=0}^{\infty}
C_{j_k\ldots j_1}\biggl|_{j_{g_1}=j_{g_2},\ldots, j_{g_{2r-1}}=j_{g_{2r}}}
$$

\noindent
does not depend on the method of summation$,$ i.e.$,$ for example$,$ for $r=3$}
$$
\sum\limits_{j_{g_1}, j_{g_3},j_{g_5}=0}^{\infty}
C_{j_k\ldots j_1}\biggl|_{j_{g_1}=j_{g_2},j_{g_3}=j_{g_4},j_{g_5}=j_{g_6}}=
\sum\limits_{j_{g_1}=0}^{\infty}
\sum\limits_{j_{g_3},j_{g_5}=0}^{\infty}
C_{j_k\ldots j_1}\biggl|_{j_{g_1}=j_{g_2},j_{g_3}=j_{g_4},j_{g_5}=j_{g_6}}=
$$
$$
=
\hspace{-1mm}\sum\limits_{j_{g_1},j_{g_3}=0}^{\infty}\sum\limits_{j_{g_5}=0}^{\infty}
C_{j_k\ldots j_1}\biggl|_{j_{g_1}=j_{g_2},j_{g_{3}}=j_{g_{4}},j_{g_5}=j_{g_6}}
=
\sum\limits_{j_{g_1}=0}^{\infty}\sum\limits_{j_{g_3}=0}^{\infty}\sum\limits_{j_{g_{5}}=0}^{\infty}
C_{j_k\ldots j_1}\biggl|_{j_{g_1}=j_{g_2},j_{g_{3}}=j_{g_{4}},j_{g_5}=j_{g_6}}\hspace{-0.5mm}.
$$

\vspace{2mm}

Suppose that the condition $(\star)$ is fulfilled.
Then,
by iteratively applying the formula (\ref{2026may10sss}), we obtain the following
inverse version of the Hu--Meyer
formula (\ref{2026may10})
$$
J'[\Phi]_{T,t}^{(i_1\ldots i_k)}
=
\hat J^{S}[\Phi]_{T,t}^{(i_1\ldots i_k)}+
$$

\vspace{-3mm}
$$
+
\sum\limits_{r=1}^{[k/2]}
(-1)^r\sum_{\stackrel{(\{\{g_1, g_2\}, \ldots, 
\{g_{2r-1}, g_{2r}\}\}, \{q_1, \ldots, q_{k-2r}\})}
{{}_{\{g_1, g_2, \ldots, 
g_{2r-1}, g_{2r}, q_1, \ldots, q_{k-2r}\}=\{1, 2, \ldots, k\}}}}
\prod\limits_{s=1}^r
{\bf 1}_{\{i_{g_{{}_{2s-1}}}=~i_{g_{{}_{2s}}}\ne 0\}}\times
$$

\vspace{4mm}
\begin{equation}
\label{2026may38abcd}
\times
\hat J^{S}\hspace{-1mm}\left[\tilde T^{k-2r}_{g_1,g_2,\ldots,g_{2r-1}, 
g_{2r}}\Phi\right]_{T,t}^{(i_{q_1}\ldots i_{q_{k-2r}})}\ \ \ \hbox{w.~p.~1.}
\end{equation}

\vspace{3mm}

It is interesting to compare the formulas 
(\ref{2026may10}) and (\ref{2026may38abcd}) with similar
formulas (\ref{30.4}) and (\ref{2024str11}) that connect
the iterated Stratonovich and It\^{o} stochastic integrals.

Suppose that all limiting traces
$T^{k-2r}_{g_1,g_2,\ldots,g_{2r-1}, g_{2r}}\Phi (t_{q_1},\ldots,t_{q_{k-2r}})$
exist
for all $r=1,2,\ldots,$ $[k/2]$ and for all possible
$g_1,g_2,$ $\ldots,g_{2r-1},g_{2r}$ (see {\rm (\ref{leto5007after})).
Recall the equality $T^{k-2r}_{g_1,g_2,\ldots,g_{2r-1}, g_{2r}}\Phi (t_{q_1},\ldots,t_{q_{k-2r}})=
\tilde T^{k-2r}_{g_1,g_2,\ldots,g_{2r-1}, g_{2r}}\Phi (t_{q_1},\ldots,t_{q_{k-2r}})$ (in 
$L_2([t, T]^{k-2r})$), which is true in this case. 
Then
\begin{equation}
\label{2026may33}
\lim\limits_{p\to\infty}
\left\Vert T^{k-2r,\hspace{0.2mm}p}_{g_1,g_2,\ldots,g_{2r-1}, 
g_{2r}}\Phi-
\tilde T^{k-2r,\hspace{0.2mm}p}_{g_1,g_2,\ldots,g_{2r-1}, 
g_{2r}}\Phi\right\Vert_{L_2([t, T]^{k-2r})}=0
\end{equation}

\noindent
for all $r=1,2,\ldots,[k/2]$ and for all possible
$g_1,g_2,\ldots,g_{2r-1},g_{2r}$.

Note that

\newpage
\noindent
$$
\left\Vert T^{k-2r,\hspace{0.2mm}p}_{g_1,g_2,\ldots,g_{2r-1}, 
g_{2r}}\Phi-
\tilde T^{k-2r,\hspace{0.2mm}p}_{g_1,g_2,\ldots,g_{2r-1}, 
g_{2r}}\Phi\right\Vert_{L_2([t, T]^{k-2r})}^2=
$$
\begin{equation}
\label{2026may32}
=\sum\limits_{j_{q_1},\ldots,j_{q_{k-2r}}=0}^p
\Biggl(
\sum\limits_{j_{g_1}, j_{g_3},\ldots ,j_{g_{2r-1}}=0}^p
C_{j_k\ldots j_1}\biggl|_{j_{g_1}=j_{g_2},\ldots, j_{g_{2r-1}}=j_{g_{2r}}}-
\tilde C_{j_{q_{k-2r}}\ldots j_{q_1}}
\Biggr)^2,
\end{equation}

\noindent
where $\tilde C_{j_{q_{k-2r}}\ldots j_{q_1}}$ is defined by (\ref{october202621}).

Now let us assume that 
all limiting traces
$\tilde T^{k-2r}_{g_1,g_2,\ldots,g_{2r-1}, g_{2r}}\Phi (t_{q_1},\ldots,t_{q_{k-2r}})$
exist
for all $r=1,2,\ldots,$ $[k/2]$ and for all possible
$g_1,g_2,$ $\ldots,g_{2r-1},g_{2r}$ (see {\rm (\ref{leto5007after})), and,
in addition, that the equality (\ref{2026may33}) is satisfied.

Then, by virtue of the equality
$$
T^{k-2r,\hspace{0.2mm}p}_{g_1,g_2,\ldots,g_{2r-1}, 
g_{2r}}\Phi=
$$

\vspace{-2mm}
$$
=\left(T^{k-2r,\hspace{0.2mm}p}_{g_1,g_2,\ldots,g_{2r-1}, 
g_{2r}}\Phi-
\tilde T^{k-2r,\hspace{0.2mm}p}_{g_1,g_2,\ldots,g_{2r-1}, 
g_{2r}}\Phi\right)+\tilde T^{k-2r,\hspace{0.2mm}p}_{g_1,g_2,\ldots,g_{2r-1}, 
g_{2r}}\Phi,
$$

\vspace{2mm}
\noindent
we obtain that 
all limiting traces
$T^{k-2r}_{g_1,g_2,\ldots,g_{2r-1}, g_{2r}}\Phi (t_{q_1},\ldots,t_{q_{k-2r}})$
exist
for all $r=1,2,\ldots,$ $[k/2]$ and for all possible
$g_1,g_2,$ $\ldots,g_{2r-1},g_{2r},$ and 
the equality $T^{k-2r}_{g_1,g_2,\ldots,g_{2r-1}, g_{2r}}\Phi (t_{q_1},\ldots,t_{q_{k-2r}})=
\tilde T^{k-2r}_{g_1,g_2,\ldots,g_{2r-1}, g_{2r}}\Phi (t_{q_1},\ldots,t_{q_{k-2r}})$ 
is fulfilled (in 
$L_2([t, T]^{k-2r})$).

Further, we will consider the connection of the above results
with Theorems~4 and 7. Suppose that the condition (\ref{2026may33})
is fulfilled
and all limiting traces 
$\tilde T^{k-2r}_{g_1,g_2,\ldots,g_{2r-1}, g_{2r}}\Phi(t_{q_1},\ldots,t_{q_{k-2r}})$ exist 
for all $r=1,2,\ldots,[k/2]$ and for all possible
$g_1,g_2,\ldots,g_{2r-1},g_{2r}.$

Using (\ref{2026may2}), we get
$$
{\sf M}\left\{\left(J'\hspace{-1mm}\left[T^{k-2r,\hspace{0.2mm}p}_{g_1,g_2,\ldots,g_{2r-1}, 
g_{2r}}\Phi\right]_{T,t}^{(i_{q_1}\ldots i_{q_{k-2r}})}-
J'\hspace{-1mm}\left[\tilde T^{k-2r,\hspace{0.2mm}p}_{g_1,g_2,\ldots,g_{2r-1}, 
g_{2r}}\Phi\right]_{T,t}^{(i_{q_1}\ldots i_{q_{k-2r}})}\right)^2\right\}\le
$$
$$
\le K_{k-2r} \left\Vert T^{k-2r,\hspace{0.2mm}p}_{g_1,g_2,\ldots,g_{2r-1}, 
g_{2r}}\Phi-
\tilde T^{k-2r,\hspace{0.2mm}p}_{g_1,g_2,\ldots,g_{2r-1}, 
g_{2r}}\Phi\right\Vert_{L_2([t, T]^{k-2r})}^2\ \to\ 0,
$$

\vspace{-2mm}
$$
{\sf M}\left\{\left(J'\hspace{-1mm}\left[\tilde T^{k-2r}_{g_1,g_2,\ldots,g_{2r-1}, 
g_{2r}}\Phi\right]_{T,t}^{(i_{q_1}\ldots i_{q_{k-2r}})}-
J'\hspace{-1mm}\left[\tilde T^{k-2r,\hspace{0.2mm}p}_{g_1,g_2,\ldots,g_{2r-1}, 
g_{2r}}\Phi\right]_{T,t}^{(i_{q_1}\ldots i_{q_{k-2r}})}\right)^2\right\}\le
$$
$$
\le K'_{k-2r} \left\Vert \tilde T^{k-2r}_{g_1,g_2,\ldots,g_{2r-1}, 
g_{2r}}\Phi-
\tilde T^{k-2r,\hspace{0.2mm}p}_{g_1,g_2,\ldots,g_{2r-1}, 
g_{2r}}\Phi\right\Vert_{L_2([t, T]^{k-2r})}^2\ \to\ 0
$$

\vspace{2mm}
\noindent
if $p\to\infty$ (for all $r=1,2,\ldots,[k/2]$ and for all possible
$g_1,g_2,\ldots,g_{2r-1},g_{2r}$
(see (\ref{leto5007after}))), where $K_{k-2r}, K'_{k-2r}$ are constants.
Then w.~p.~1
$$
\hbox{\vtop{\offinterlineskip\halign{
\hfil#\hfil\cr
{\rm l.i.m.}\cr
$\stackrel{}{{}_{p\to \infty}}$\cr
}} }
J'\hspace{-1mm}\left[T^{k-2r,\hspace{0.2mm}p}_{g_1,g_2,\ldots,g_{2r-1}, 
g_{2r}}\Phi\right]_{T,t}^{(i_{q_1}\ldots i_{q_{k-2r}})}=
$$

\vspace{-2mm}
$$
=
\hbox{\vtop{\offinterlineskip\halign{
\hfil#\hfil\cr
{\rm l.i.m.}\cr
$\stackrel{}{{}_{p\to \infty}}$\cr
}} }
J'\hspace{-1mm}\left[\tilde T^{k-2r,\hspace{0.2mm}p}_{g_1,g_2,\ldots,g_{2r-1}, 
g_{2r}}\Phi\right]_{T,t}^{(i_{q_1}\ldots i_{q_{k-2r}})}=
$$

\vspace{-2mm}
\begin{equation}
\label{october202615}
=
J'\hspace{-1mm}\left[\tilde T^{k-2r}_{g_1,g_2,\ldots,g_{2r-1}, 
g_{2r}}\Phi\right]_{T,t}^{(i_{q_1}\ldots i_{q_{k-2r}})}.
\end{equation}

\vspace{2mm}

Passing to the limit
$\hbox{\vtop{\offinterlineskip\halign{
\hfil#\hfil\cr
{\rm l.i.m.}\cr
$\stackrel{}{{}_{p\to \infty}}$\cr
}} }$ in (\ref{2026may1}) and applying (\ref{october202615}), we obtain
w.~p.~1

\vspace{-2mm}
$$
\hat J^{S}[\Phi]_{T,t}^{(i_1\ldots i_k)}=\hbox{\vtop{\offinterlineskip\halign{
\hfil#\hfil\cr
{\rm l.i.m.}\cr
$\stackrel{}{{}_{p\to \infty}}$\cr
}} }\hat J^{S}_p[\Phi]_{T,t}^{(i_1\ldots i_k)}
=
\hbox{\vtop{\offinterlineskip\halign{
\hfil#\hfil\cr
{\rm l.i.m.}\cr
$\stackrel{}{{}_{p\to \infty}}$\cr
}} }J'[\Phi_p]_{T,t}^{(i_1\ldots i_k)}+
$$

\vspace{-3mm}
$$
+
\sum\limits_{r=1}^{[k/2]}
\sum_{\stackrel{(\{\{g_1, g_2\}, \ldots, 
\{g_{2r-1}, g_{2r}\}\}, \{q_1, \ldots, q_{k-2r}\})}
{{}_{\{g_1, g_2, \ldots, 
g_{2r-1}, g_{2r}, q_1, \ldots, q_{k-2r}\}=\{1, 2, \ldots, k\}}}}
\prod\limits_{s=1}^r
{\bf 1}_{\{i_{g_{{}_{2s-1}}}=~i_{g_{{}_{2s}}}\ne 0\}}\times
$$

\vspace{4mm}
\begin{equation}
\label{october202616}
\times
\hbox{\vtop{\offinterlineskip\halign{
\hfil#\hfil\cr
{\rm l.i.m.}\cr
$\stackrel{}{{}_{p\to \infty}}$\cr
}} }J'\hspace{-1mm}\left[\tilde T^{k-2r,\hspace{0.2mm}p}_{g_1,g_2,\ldots,g_{2r-1}, 
g_{2r}}\Phi\right]_{T,t}^{(i_{q_1}\ldots i_{q_{k-2r}})},
\end{equation}

\vspace{4mm}
\noindent
where notations are the same as in (\ref{2026may1}).

Replacing the function $\Phi(t_1,\ldots,t_k)$ in (\ref{october202616}) 
with the Volterra-type kernel (see (\ref{july7000}))
\begin{equation}
\label{october202633}
K(t_1,\ldots,t_k)=
\left\{\begin{matrix}
\psi_1(t_1)\ldots \psi_k(t_k),\ &t_1<\ldots<t_k\cr\cr
0,\ &\hbox{\rm otherwise}
\end{matrix}
\right.\ \ \ (k\ge 2),
\end{equation}

\noindent
where $\psi_1(\tau),\ldots,\psi_k(\tau)\in L_2([t,T])$ and
$t_1,\ldots,t_k\in [t, T],$ and using (\ref{febr5000}), we have

\newpage
\noindent
$$
\hbox{\vtop{\offinterlineskip\halign{
\hfil#\hfil\cr
{\rm l.i.m.}\cr
$\stackrel{}{{}_{p\to \infty}}$\cr
}} }\sum_{j_1,\ldots,j_k=0}^{p}
C_{j_k\ldots j_1}
\prod_{l=1}^k \zeta_{j_l}^{(i_l)}
=
J[\psi^{(k)}]_{T,t}^{(i_1\ldots i_k)}+
$$

\vspace{-2mm}
$$
+
\sum\limits_{r=1}^{[k/2]}
\sum_{\stackrel{(\{\{g_1, g_2\}, \ldots, 
\{g_{2r-1}, g_{2r}\}\}, \{q_1, \ldots, q_{k-2r}\})}
{{}_{\{g_1, g_2, \ldots, 
g_{2r-1}, g_{2r}, q_1, \ldots, q_{k-2r}\}=\{1, 2, \ldots, k\}}}}
\prod\limits_{s=1}^r
{\bf 1}_{\{i_{g_{{}_{2s-1}}}=~i_{g_{{}_{2s}}}\ne 0\}}\times
$$

\vspace{5mm}
\begin{equation}
\label{october202617}
\times
\hbox{\vtop{\offinterlineskip\halign{
\hfil#\hfil\cr
{\rm l.i.m.}\cr
$\stackrel{}{{}_{p\to \infty}}$\cr
}} }J'\hspace{-1mm}\left[\tilde T^{k-2r,\hspace{0.2mm}p}_{g_1,g_2,\ldots,g_{2r-1}, 
g_{2r}}K\right]_{T,t}^{(i_{q_1}\ldots i_{q_{k-2r}})},
\end{equation}

\vspace{4mm}
\noindent
where $J[\psi^{(k)}]_{T,t}^{(i_1\ldots i_k)}$
is the iterated It\^{o} stochastic integral
(\ref{ito}).

Consider the equality (see the proof of Theorem~2.49 in \cite{2018a}, Sect.~2.22)
$$
\hbox{\vtop{\offinterlineskip\halign{
\hfil#\hfil\cr
{\rm l.i.m.}\cr
$\stackrel{}{{}_{p\to \infty}}$\cr
}} }
\sum\limits_{j_{q_1},\ldots,j_{q_{k-2r}}=0}^p
\frac{1}{2^r}
C_{j_k \ldots j_1}\biggl|_{(j_{g_2} j_{g_1})\curvearrowright (\cdot)
\ldots (j_{g_{2r}} j_{g_{2r-1}})\curvearrowright (\cdot),
j_{g_{{}_{1}}}=~j_{g_{{}_{2}}},\ldots, j_{g_{{}_{2r-1}}}=~j_{g_{{}_{2r}}}}\biggr.
\times 
$$

\vspace{2mm}
$$
\times
\prod\limits_{s=1}^r
{\bf 1}_{\{i_{g_{{}_{2s-1}}}=~i_{g_{{}_{2s}}}\ne 0\}}
J'[\phi_{j_{q_1}}\ldots \phi_{j_{q_{k-2r}}}]_{T,t}^{(i_{q_1}\ldots i_{q_{k-2r}})}=
$$

\vspace{1mm}
\begin{equation}
\label{october202618}
=\frac{1}{2^r}
J[\psi^{(k)}]_{T,t}^{s_r, \ldots, s_1}
\end{equation}

\vspace{1mm}
\noindent
w.~p.~1, where $g_{2}=g_{1}+1,\ldots, g_{2r}=g_{2r-1}+1,$
$g_{2i-1}\stackrel{\sf def}{=}s_i;$\ $i=1,2,\ldots,r;$\
$r=1,2,\ldots,\left[k/2\right],$ 
$(s_r,\ldots,s_1)\in {\rm A}_{k,r},$ $J[\psi^{(k)}]_{T,t}^{s_r,\ldots,s_1}$ is
defined by (\ref{30.1}) and ${\rm A}_{k,r}$ is defined by (\ref{30.5550001});
another notations in (\ref{october202618}) are the same as in Sect.~2.1.

Also, consider the equality (see the proof in \cite{2018a}, Sect.~2.27.4 and 2.30)
$$
\lim\limits_{p\to\infty}
\sum\limits_{j_{g_1}, j_{g_3},\ldots ,j_{g_{2r-1}}=0}^p
C_{j_k\ldots j_1}\biggl|_{j_{g_1}=j_{g_2},\ldots, j_{g_{2r-1}}=j_{g_{2r}}}=
$$
\begin{equation}
\label{october202619}
=\frac{1}{2^r} \prod\limits_{l=1}^r {\bf 1}_{\{g_{2l}=g_{2l-1}+1\}}
C_{j_k \ldots j_1}\biggl|_{(j_{g_2} j_{g_1})\curvearrowright (\cdot)
\ldots (j_{g_{2r}} j_{g_{2r-1}})\curvearrowright (\cdot),
j_{g_{{}_{1}}}=~j_{g_{{}_{2}}},\ldots, j_{g_{{}_{2r-1}}}=~j_{g_{{}_{2r}}}
}\biggr.
\end{equation}

\vspace{3mm}
\noindent
for all possible $g_1,g_2,\ldots,g_{2r-1},g_{2r}$ (see {\rm (\ref{leto5007after})),
where $k\ge 2r,$ $r=1,2,\ldots,$ $[k/2],$ $C_{j_k\ldots j_1}$ is defined by (\ref{after3000});
another notations are the same as in Sect.~2.1. 

Note that
$$
\sum_{\stackrel{(\{\{g_1, g_2\}, \ldots, 
\{g_{2r-1}, g_{2r}\}\}, \{q_1, \ldots, q_{k-2r}\})}
{{}_{\{g_1, g_2, \ldots, 
g_{2r-1}, g_{2r}, q_1, \ldots, q_{k-2r}\}=\{1, 2, \ldots, k\}}}}
\Biggl|_{g_2=g_1+1, g_3=g_2+1,\ldots, g_{2r}=g_{2r-1}+1}\Biggr.
A_{g_1,g_3,\ldots,g_{2r-1}}=
$$

\vspace{-3mm}
\begin{equation}
\label{after800}
=\sum\limits_{(s_r,\ldots,s_1)\in {\rm A}_{k,r}}
A_{s_1,s_2,\ldots,s_r},
\end{equation}

\vspace{2mm}
\noindent
where $A_{g_1,g_3,\ldots,g_{2r-1}},$ 
$A_{s_1,s_2,\ldots,s_r}$ are scalar values,
$g_{2i-1}=s_i;$\ $i=1,2,\ldots,r;$\ $r=1,2,\ldots,\left[k/2\right],$
${\rm A}_{k,r}$ is defined by (\ref{30.5550001}).

Combining (\ref{october202620}), (\ref{october202621}), (\ref{october202618}),
(\ref{october202619}) and (\ref{after800}), we obtain w.~p.~1

\vspace{-3mm}
$$
\prod\limits_{s=1}^r
{\bf 1}_{\{i_{g_{{}_{2s-1}}}=~i_{g_{{}_{2s}}}\ne 0\}}\
\hbox{\vtop{\offinterlineskip\halign{
\hfil#\hfil\cr
{\rm l.i.m.}\cr
$\stackrel{}{{}_{p\to \infty}}$\cr
}} }J'\hspace{-1mm}\left[\tilde T^{k-2r,\hspace{0.2mm}p}_{g_1,g_2,\ldots,g_{2r-1}, 
g_{2r}}K\right]_{T,t}^{(i_{q_1}\ldots i_{q_{k-2r}})}=
$$
$$
=\prod\limits_{l=1}^r {\bf 1}_{\{g_{2l}=g_{2l-1}+1\}}
\frac{1}{2^r}
J[\psi^{(k)}]_{T,t}^{s_r, \ldots, s_1},
$$

\noindent
and
$$
J[\psi^{(k)}]_{T,t}^{(i_1\ldots i_k)}+
$$

\vspace{-4mm}
$$
+\sum\limits_{r=1}^{[k/2]}
\sum_{\stackrel{(\{\{g_1, g_2\}, \ldots, 
\{g_{2r-1}, g_{2r}\}\}, \{q_1, \ldots, q_{k-2r}\})}
{{}_{\{g_1, g_2, \ldots, 
g_{2r-1}, g_{2r}, q_1, \ldots, q_{k-2r}\}=\{1, 2, \ldots, k\}}}}
\prod\limits_{s=1}^r
{\bf 1}_{\{i_{g_{{}_{2s-1}}}=~i_{g_{{}_{2s}}}\ne 0\}}\times
$$

\vspace{2mm}
$$
\times
\hbox{\vtop{\offinterlineskip\halign{
\hfil#\hfil\cr
{\rm l.i.m.}\cr
$\stackrel{}{{}_{p\to \infty}}$\cr
}} }J'\hspace{-1mm}\left[\tilde T^{k-2r,\hspace{0.2mm}p}_{g_1,g_2,\ldots,g_{2r-1}, 
g_{2r}}K\right]_{T,t}^{(i_{q_1}\ldots i_{q_{k-2r}})}=
$$

\begin{equation}
\label{october202623}
=J[\psi^{(k)}]_{T,t}^{(i_1\ldots i_k)}+
\sum_{r=1}^{\left[k/2\right]}\frac{1}{2^r}
\sum_{(s_r,\ldots,s_1)\in {\rm A}_{k,r}}
J[\psi^{(k)}]_{T,t}^{s_r,\ldots,s_1}.
\end{equation}

\vspace{2mm}

Suppose that 
$\psi_1(\tau), \ldots, \psi_k(\tau)$ are continuous 
functions on $[t, T].$ Further, applying Theorem~5 
to the right-hand side of (\ref{october202623})
and combining (\ref{october202617}), (\ref{october202623}), we get
the following expansion
$$
J^{*}[\psi^{(k)}]_{T,t}^{(i_1\ldots i_k)}=\hbox{\vtop{\offinterlineskip\halign{
\hfil#\hfil\cr
{\rm l.i.m.}\cr
$\stackrel{}{{}_{p\to \infty}}$\cr
}} }\sum_{j_1,\ldots,j_k=0}^{p}
C_{j_k\ldots j_1}
\prod_{l=1}^k \zeta_{j_l}^{(i_l)},
$$

\noindent
where $J^{*}[\psi^{(k)}]_{T,t}^{(i_1\ldots i_k)}$ is the iterated
Stratonovich stochastic integral (\ref{str}).

Note that in our case (the case of Volterra-type kernel (\ref{october202633})), limiting traces
$\tilde T^{k-2r}_{g_1,g_2,\ldots,g_{2r-1}, g_{2r}}K$ exist 
for all $r=1,2,\ldots,[k/2]$ and for all possible
$g_1,g_2,\ldots,g_{2r-1},g_{2r},$ and are determined
by the equality (\ref{2025may25}) (also see (\ref{2025may28})).
Moreover, the condition (\ref{2026may33}) in this case
has the form (\ref{july700000}) (or (\ref{2025may26})).

Further, we will consider another approach to deriving
the Hu--Meyer formulas based on the so-called
Hilbert space valued traces \cite{HuHu}, \cite{bugh1}, \cite{bugh3}. At that we will still
consider the case of a multidimensional Wiener process.

Recall the following definition of multiple Stratonovich stochastic
integral (see, for example, \cite{bugh1})
\begin{equation}
\label{30.34ququsss}
\hbox{\vtop{\offinterlineskip\halign{
\hfil#\hfil\cr
{\rm l.i.m.}\cr
$\stackrel{}{{}_{N\to \infty}}$\cr
}} }\sum_{j_1,\ldots,j_k=0}^{N-1}
\Phi\left(\tau_{j_1},\ldots,\tau_{j_k}\right)
\prod\limits_{l=1}^k\Delta{\bf w}_{\tau_{j_l}}^{(i_l)}
\stackrel{\rm def}{=}J[\Phi]_{T,t}^{(i_1\ldots i_k)},
\end{equation}

\noindent
where 
$\Phi(t_1,\ldots,t_k):\ [t, T]^k\to{\bf R}$ is a 
continuous nonrandom
function, 
${\bf w}_{\tau}$ is a random vector with 
an $m+1$ components
(${\bf w}_{\tau}^{(i)} $ $(i=1,\ldots,m)$
are independent standard Wiener processes and
${\bf w}_{\tau}^{(0)}=\tau$),
$\Delta {\bf w}_{\tau_j}^{(i)}={\bf w}_{\tau_{j+1}}^{(i)}-{\bf w}_{\tau_j}^{(i)},$
$\left\{\tau_{j}\right\}_{j=0}^{N}$ is a partition of
$[t,T]$ which satisfies the following condition 
\begin{equation}
\label{novem2026xxx6}
t=\tau_0<\tau_1<\ldots <\tau_N=T,\ \ \
\max\limits_{0\le j\le N-1}\left|\tau_{j+1}-\tau_j\right|\to 0\ \
\hbox{if}\ \ N\to \infty.
\end{equation}

Note that the function $\Phi(t_1,\ldots,t_k)$
in (\ref{30.34ququsss}) can be discontinuous.
For example, it can have the form (2.379) (see \cite{2018a}, Sect.~2.4.1).

Further, let us write the formula (2.429) (see \cite{2018a}, Sect.~2.4.1), replacing  
$R_{p_1 p_2 p_3 p_4}(t_1,\ldots,t_4)$ with
$\Phi(t_1,\ldots,t_4)$ and using the formula (1.329) (see \cite{2018a}, Sect.~1.11) for the
multiple Wiener stochastic integral $J'[\Phi]_{T,t}^{(i_1\ldots i_4)}$ defined by (1.328)
(see \cite{2018a}, Sect.~1.11)

\vspace{-2mm}
$$
J[\Phi]_{T,t}^{(i_1\ldots i_4)}=J'[\Phi]_{T,t}^{(i_1\ldots i_4)}+
$$

\vspace{-3mm}
$$
+{\bf 1}_{\{i_1=i_2\ne 0\}}J'\left[\Phi(t_1,\ldots,t_4)\bigl.\bigr|_{t_1=t_2}
\right]_{T,t}^{(0 i_3 i_4)}+
{\bf 1}_{\{i_1=i_3\ne 0\}}J'\left[\Phi(t_1,\ldots,t_4)\bigl.\bigr|_{t_1=t_3}
\right]_{T,t}^{(0 i_2 i_4)}+
$$
$$
+{\bf 1}_{\{i_1=i_4\ne 0\}}J'\left[\Phi(t_1,\ldots,t_4)\bigl.\bigr|_{t_1=t_4}
\right]_{T,t}^{(0 i_2 i_3)}+
{\bf 1}_{\{i_2=i_3\ne 0\}}J'\left[\Phi(t_1,\ldots,t_4)\bigl.\bigr|_{t_2=t_3}
\right]_{T,t}^{(i_1 0 i_4)}+
$$

\vspace{-3mm}
$$
+{\bf 1}_{\{i_2=i_4\ne 0\}}J'\left[\Phi(t_1,\ldots,t_4)\bigl.\bigr|_{t_2=t_4}
\right]_{T,t}^{(i_1 0 i_3)}+
{\bf 1}_{\{i_3=i_4\ne 0\}}J'\left[\Phi(t_1,\ldots,t_4)\bigl.\bigr|_{t_3=t_4}
\right]_{T,t}^{(i_1 i_2 0)}+
$$

\vspace{-1mm}
$$
+{\bf 1}_{\{i_1=i_2\ne 0\}}{\bf 1}_{\{i_3=i_4\ne 0\}}
J'\left[\Phi(t_1,\ldots,t_4)\bigl.\bigr|_{t_1=t_2,t_3=t_4}
\right]_{T,t}^{(00)}+
$$

\vspace{-1mm}
$$
+{\bf 1}_{\{i_1=i_3\ne 0\}}{\bf 1}_{\{i_2=i_4\ne 0\}}
J'\left[\Phi(t_1,\ldots,t_4)\bigl.\bigr|_{t_1=t_3,t_2=t_4}
\right]_{T,t}^{(00)}+
$$

\vspace{-3mm}
\begin{equation}
\label{2025nov100}
+{\bf 1}_{\{i_1=i_4\ne 0\}}{\bf 1}_{\{i_2=i_3\ne 0\}}
J'\left[\Phi(t_1,\ldots,t_4)\bigl.\bigr|_{t_1=t_4,t_2=t_3}
\right]_{T,t}^{(00)}
\end{equation}

\vspace{2mm}
\noindent
w.~p.~1. Further, applying a Fubini-type theorem for stochastic and 
Lebesgue intgrals (as in \cite{farre}) to the right-hand side of
(\ref{2025nov100}), we obtain w.~p.~1

\vspace{-2mm}
$$
J[\Phi]_{T,t}^{(i_1\ldots i_4)}=
J'[\Phi]_{T,t}^{(i_1\ldots i_4)}
+
$$

\vspace{-2mm}
$$
+{\bf 1}_{\{i_1=i_2\ne 0\}}J'\Biggl[
~\int\limits_{[t, T]} \Phi(t_1,\ldots,t_4)\bigl.\bigr|_{t_1=t_2}
dt_1\Biggr]_{T,t}^{(i_3 i_4)}+
$$

\vspace{-2mm}
$$
+
{\bf 1}_{\{i_1=i_3\ne 0\}}J'\Biggl[~\int\limits_{[t, T]}\Phi(t_1,\ldots,t_4)\bigl.\bigr|_{t_1=t_3}
dt_1
\Biggr]_{T,t}^{(i_2 i_4)}+
$$

\vspace{-2mm}
$$
+{\bf 1}_{\{i_1=i_4\ne 0\}}J'\Biggl[~\int\limits_{[t, T]}
\Phi(t_1,\ldots,t_4)\bigl.\bigr|_{t_1=t_4}
dt_1\Biggr]_{T,t}^{(i_2 i_3)}+
$$

\vspace{-2mm}
$$
+
{\bf 1}_{\{i_2=i_3\ne 0\}}J'\Biggl[~\int\limits_{[t, T]}
\Phi(t_1,\ldots,t_4)\bigl.\bigr|_{t_2=t_3}
dt_2\Biggr]_{T,t}^{(i_1 i_4)}+
$$

\vspace{-2mm}
$$
+{\bf 1}_{\{i_2=i_4\ne 0\}}J'\Biggl[~\int\limits_{[t, T]}
\Phi(t_1,\ldots,t_4)\bigl.\bigr|_{t_2=t_4}
dt_2\Biggr]_{T,t}^{(i_1 i_3)}+
$$
$$
+
{\bf 1}_{\{i_3=i_4\ne 0\}}J'\Biggl[~\int\limits_{[t, T]}
\Phi(t_1,\ldots,t_4)\bigl.\bigr|_{t_3=t_4}
dt_3\Biggr]_{T,t}^{(i_1 i_2)}+
$$

$$
+{\bf 1}_{\{i_1=i_2\ne 0\}}{\bf 1}_{\{i_3=i_4\ne 0\}}
\int\limits_{[t, T]^2}\Phi(t_1,\ldots,t_4)\bigl.\bigr|_{t_1=t_2,t_3=t_4}
dt_1 dt_3+
$$

$$
+{\bf 1}_{\{i_1=i_3\ne 0\}}{\bf 1}_{\{i_2=i_4\ne 0\}}
\int\limits_{[t, T]^2}\Phi(t_1,\ldots,t_4)\bigl.\bigr|_{t_1=t_3,t_2=t_4}
dt_1 dt_2+
$$

\vspace{-2mm}
\begin{equation}
\label{2025nov101}
+{\bf 1}_{\{i_1=i_4\ne 0\}}{\bf 1}_{\{i_2=i_3\ne 0\}}
\int\limits_{[t, T]^2}\Phi(t_1,\ldots,t_4)\bigl.\bigr|_{t_1=t_4,t_2=t_3}
dt_1 dt_2.
\end{equation}

\vspace{2mm}

Let us generalize the formula (\ref{2025nov101}) to the case $k\in{\bf N}$.
As a result, we obtain the following Hu--Meyer formula for the 
case of a multidimensional Wiener process
$$
J[\Phi]_{T,t}^{(i_1\ldots i_k)}
=
J'[\Phi]_{T,t}^{(i_1\ldots i_k)}+
$$

\vspace{-3mm}
$$
+
\sum\limits_{r=1}^{[k/2]}
\sum_{\stackrel{(\{\{g_1, g_2\}, \ldots, 
\{g_{2r-1}, g_{2r}\}\}, \{q_1, \ldots, q_{k-2r}\})}
{{}_{\{g_1, g_2, \ldots, 
g_{2r-1}, g_{2r}, q_1, \ldots, q_{k-2r}\}=\{1, 2, \ldots, k\}}}}
\prod\limits_{s=1}^r
{\bf 1}_{\{i_{g_{{}_{2s-1}}}=~i_{g_{{}_{2s}}}\ne 0\}}\times
$$

\vspace{4mm}
\begin{equation}
\label{2026may10zyx}
\times
J'\hspace{-1mm}\left[\breve T^{k-2r}_{g_1,g_2,\ldots,g_{2r-1}, 
g_{2r}}\Phi\right]_{T,t}^{(i_{q_1}\ldots i_{q_{k-2r}})}\ \ \ \hbox{w.~p.~1,}
\end{equation}

\vspace{3mm}
\noindent
where the so-called 
Hilbert space valued trace 
$\breve T^{k-2r}_{g_1,g_2,\ldots,g_{2r-1}, 
g_{2r}}\Phi\in L_2([t, T]^{k-2r})$ 
has the form (we assume that this trace exists)

\vspace{-2mm}
$$
\breve T^{k-2r}_{g_1,g_2,\ldots,g_{2r-1}, 
g_{2r}}\Phi(t_{q_1},\ldots,t_{q_{k-2r}})\stackrel{\sf def}{=}
$$

$$
\stackrel{\sf def}{=}
\int\limits_{[t, T]^r}
\Phi(
t_1,\ldots,t_k)\biggl.\biggr|_{t_{g_{{}_{1}}}=t_{g_{{}_{2}}},\ldots,
t_{g_{{}_{2r-1}}}=t_{g_{{}_{2r}}}}dt_{g_{{}_{1}}}\ldots dt_{g_{{}_{2r-1}}},
$$

\vspace{2mm}
\noindent
where 
$$
J'\hspace{-1mm}\left[\breve T^{k-2r}_{g_1,g_2,\ldots,g_{2r-1}, 
g_{2r}}\Phi\right]_{T,t}^{(i_{q_1}\ldots i_{q_{k-2r}})}\stackrel{\sf def}{=}
\breve T^{k-2r}_{g_1,g_2,\ldots,g_{2r-1}, 
g_{2r}}\Phi
$$ 

\vspace{1mm}
\noindent
for $k=2r;$
another notations are the same as in (\ref{2026may10}).
Recall that
in \cite{bugh1}, the variant $\breve T^{k-2r}_{1,2\ldots,2r-1,2r}\Phi(t_{2r+1},\ldots,t_k)$
of $\breve T^{k-2r}_{g_1,g_2,\ldots,g_{2r-1}, 
g_{2r}}\Phi(t_{q_1},\ldots,t_{q_{k-2r}})$ was considered.

The inverse version of the Hu--Meyer
formula (\ref{2026may10zyx}) has the form

\vspace{-4mm}
$$
J'[\Phi]_{T,t}^{(i_1\ldots i_k)}
=
J[\Phi]_{T,t}^{(i_1\ldots i_k)}+
$$

\vspace{-3mm}
$$
+
\sum\limits_{r=1}^{[k/2]}
(-1)^r\sum_{\stackrel{(\{\{g_1, g_2\}, \ldots, 
\{g_{2r-1}, g_{2r}\}\}, \{q_1, \ldots, q_{k-2r}\})}
{{}_{\{g_1, g_2, \ldots, 
g_{2r-1}, g_{2r}, q_1, \ldots, q_{k-2r}\}=\{1, 2, \ldots, k\}}}}
\prod\limits_{s=1}^r
{\bf 1}_{\{i_{g_{{}_{2s-1}}}=~i_{g_{{}_{2s}}}\ne 0\}}\times
$$

\vspace{4mm}
$$
\times
J\hspace{-1mm}\left[\breve T^{k-2r}_{g_1,g_2,\ldots,g_{2r-1}, 
g_{2r}}\Phi\right]_{T,t}^{(i_{q_1}\ldots i_{q_{k-2r}})}\ \ \ \hbox{w.~p.~1.}
$$

\vspace{3mm}

Note that there are more general 
definitions of the multiple Stratonovich stochastic integral
(see, for example, (1.5.9) and Sect.~2.1 in \cite{bugh1})
that are consistent with definition (\ref{30.34ququsss}) on the class
of continuous functions. 

Further, consider one of these generalizations. 
But first, consider the following Riemann--Stieltjes integral
\begin{equation}
\label{novem2026xxx2}
\int\limits_{[t, T]^k}\Phi(t_1,\ldots,t_k)
d{\bf w}_{t_1}^{(i_1)p}\ldots d{\bf w}_{t_k}^{(i_k)p},
\end{equation}
where $p\in{\bf N}$, $i_1,\ldots,i_k=0,1,\ldots,m,$ $\Phi(t_1,\ldots,t_k):\ [t, T]^k\to{\bf R}$ is a 
nonrandom function, 
${\bf w}_{\tau}$ is a random vector with 
an $m+1$ components
(${\bf w}_{\tau}^{(i)} $ $(i=1,\ldots,m)$
are independent standard Wiener processes and
${\bf w}_{\tau}^{(0)}=\tau$),
${\bf w}_{\tau}^{(i)p}$ is the following 
mean-square approximation of 
${\bf w}_{\tau}^{(i)}$ \cite{Lipt}
\begin{equation}
\label{novem2026xxx1000}
{\bf w}_{\tau}^{(i)p}={\bf w}_{t}^{(i)p}+
\sum_{j=0}^{p}\int\limits_t^{\tau}
\phi_j(s)ds\ \zeta_j^{(i)},
\end{equation}
where
$$
\zeta_j^{(i)}=
\int\limits_t^T \phi_j(s)d{\bf w}_s^{(i)},
$$

\noindent
$\tau\in[t, T],$ $t\ge 0,$
$\{\phi_j(x)\}_{j=0}^{\infty}$ is an arbitrary CONS 
of functions in the space $L_2([t, T]),$ and
$\zeta_j^{(i)}$ are independent standard Gaussian 
random variables for various $i$ or $j$ (in the case when $i\ne 0$).

From (\ref{novem2026xxx1000}) we have
\begin{equation}
\label{novem2026xxx1}
d{\bf w}_{\tau}^{(i)p}=
\sum_{j=0}^{p}
\phi_j(\tau)\zeta_j^{(i)} d\tau.
\end{equation}

\vspace{1mm}

Let us substitute (\ref{novem2026xxx1}) into (\ref{novem2026xxx2})
\begin{equation}
\label{novem2026xxx3}
\int\limits_{[t, T]^k}\Phi(t_1,\ldots,t_k)
d{\bf w}_{t_1}^{(i_1)p}\ldots d{\bf w}_{t_k}^{(i_k)p}
=\sum\limits_{j_1,\ldots,j_k=0}^{p}
C_{j_k \ldots j_1}\prod\limits_{l=1}^k \zeta_{j_l}^{(i_l)},
\end{equation}

\noindent
where
$$
C_{j_k \ldots j_1}=\int\limits_{[t,T]^k}
\Phi(t_1,\ldots,t_k)\prod\limits_{l=1}^k \phi_{j_l}(t_l)
dt_1\ldots dt_k
$$

\vspace{1mm}
\noindent
is the Fourier coefficient.

The formula (\ref{novem2026xxx3}) explains why
$\hat J^{S}[\Phi]_{T,t}^{(i_1\ldots i_k)}$ (see
(\ref{novem2026xxx4}), (\ref{novem2026xxx5}))
was called
the multiple Stratonovich stochastic integral.

Now consider another approximation of the Wiener process,
namely the so-called polygonal approximation (see \cite{HuHu}, Example~5.17)

\vspace{-3mm}
\begin{equation}
\label{novem2026xxx7}
{\bf w}_{\tau}^{(i)N}=
\sum\limits_{j=0}^{N-1}\left({\bf w}_{\tau_{j}}^{(i)}+
\frac{1}{\Delta_j}\Delta {\bf w}_{\tau_j}^{(i)}
(\tau-\tau_j)\right)
{\bf 1}_{T_j}(\tau),
\end{equation}

\vspace{1mm}
\noindent
where $N\in{\bf N},$ $i=0,1,\ldots,m,$  $\tau\in [t, T]$, $t\ge 0,$
${\bf w}_{\tau}$ is a random vector with 
an $m+1$ components (${\bf w}_{\tau}^{(i)} $ $(i=1,\ldots,m)$
are independent standard Wiener processes and
${\bf w}_{\tau}^{(0)}=\tau$),
$\Delta {\bf w}_{\tau_j}^{(i)}={\bf w}_{\tau_{j+1}}^{(i)}-{\bf w}_{\tau_j}^{(i)},$
$\Delta_j=\tau_{j+1}-\tau_j,$
$T_j=[\tau_j, \tau_{j+1}),$ ${\bf 1}_A$ is the indicator of the set $A,$
$\{\tau_j\}_{j=0}^N$ is a partition of
the interval $[t,T]$ which satisfies the condition (\ref{novem2026xxx6}).

From (\ref{novem2026xxx7}) we have
\begin{equation}
\label{novem2026xxx8}
d{\bf w}_{\tau}^{(i)N}=
\sum\limits_{j=0}^{N-1}
\frac{1}{\Delta_j}\Delta {\bf w}_{\tau_j}^{(i)}
d\tau
{\bf 1}_{T_j}(\tau).
\end{equation}

Using (\ref{novem2026xxx8}), we obtain
$$
\int\limits_{[t, T]^k}\Phi(t_1,\ldots,t_k)
d{\bf w}_{t_1}^{(i_1)N}\ldots d{\bf w}_{t_k}^{(i_k)N}
=\sum\limits_{j_1,\ldots,j_k=0}^{N-1}
\frac{1}{\Delta_{j_1}\ldots \Delta_{j_k}}\times
$$
$$
\times
\int\limits_{[t,T]^k}\Phi(t_1,\ldots,t_k)
{\bf 1}_{T_{j_1}}(t_1)\ldots {\bf 1}_{T_{j_k}}(t_k)dt_1\ldots dt_k
\Delta {\bf w}_{\tau_{j_1}}^{(i_1)}\ldots \Delta {\bf w}_{\tau_{j_k}}^{(i_k)}=
$$
$$
=\sum\limits_{j_1,\ldots,j_k=0}^{N-1}
\frac{1}{\Delta_{j_1}\ldots \Delta_{j_k}}
\int\limits_{T_{j_1}\times \ldots \times T_{j_k}}\Phi(t_1,\ldots,t_k) dt_1\ldots dt_k
\Delta {\bf w}_{\tau_{j_1}}^{(i_1)}\ldots \Delta {\bf w}_{\tau_{j_k}}^{(i_k)}.
$$

\vspace{2mm}

We define the multiple Stratonovich
stochastic integral (see \cite{HuHu}, \cite{SU11})
for $\Phi(t_1,\ldots,t_k)\in L_2([t, T]^k)$ (but this function 
satisfies an additional condition (see below))
as the following mean-square limit (the case of a multidimensional Wiener process)
$$
\bar J^{S}[\Phi]_{T,t}^{(i_1\ldots i_k)}=
\hbox{\vtop{\offinterlineskip\halign{
\hfil#\hfil\cr
{\rm l.i.m.}\cr
$\stackrel{}{{}_{N\to \infty}}$\cr
}} }
\sum\limits_{j_1,\ldots,j_k=0}^{N-1}
\frac{1}{\Delta_{j_1}\ldots \Delta_{j_k}}\times
$$

\begin{equation}
\label{novem2026xxx10}
\times
\int\limits_{T_{j_1}\times \ldots \times T_{j_k}}\Phi(t_1,\ldots,t_k) dt_1\ldots dt_k
\Delta {\bf w}_{\tau_{j_1}}^{(i_1)}\ldots \Delta {\bf w}_{\tau_{j_k}}^{(i_k)},
\end{equation}

\vspace{2mm}
\noindent
where $i_1,\ldots,i_k=0, 1,\ldots,m$ and $\{\tau_j\}_{j=0}^N$ 
as in (\ref{novem2026xxx6}).

We define the 
$r$th trace of the function $\Phi(t_1,\ldots,t_k)\in L_2([t, T]^k)$
by the following expression \cite{HuHu}, \cite{SU11}

\vspace{-2mm}
$$
\bar T^{k-2r}_{g_1,g_2,\ldots,g_{2r-1}, 
g_{2r}}\Phi(t_{q_1},\ldots,t_{q_{k-2r}})\stackrel{\sf def}{=}
$$

\vspace{-1mm}
\begin{equation}
\label{novem2026xxx11}
\stackrel{\sf def}{=}
\lim\limits_{N\to\infty}\bar T^{k-2r, N}_{g_1,g_2,\ldots,g_{2r-1}, 
g_{2r}}\Phi(t_{q_1},\ldots,t_{q_{k-2r}})
\end{equation}

\newpage
\noindent
in $L_2([t, T]^{k-2r})$, where 
$$
\bar T^{k-2r, N}_{g_1,g_2,\ldots,g_{2r-1}, 
g_{2r}}\Phi(t_{q_1},\ldots,t_{q_{k-2r}})=
$$
$$
=
\sum\limits_{j_{g_1},j_{g_3},\ldots,j_{g_{2r-1}}=0}^{N-1}
\frac{1}{\Delta_{j_{g_1}}\Delta_{j_{g_3}},\ldots \Delta_{j_{g_{2r-1}}}}
\times
$$

$$
\times
\int\limits_{T_{j_{g_1}}^2\times T_{j_{g_3}}^2\times
\ldots \times T_{j_{g_{2r-1}}}^2}\Phi(t_1,\ldots,t_k) dt_{g_1}dt_{g_2}\ldots dt_{g_{2r-1}}dt_{g_{2r}},
$$

\vspace{2mm}
\noindent
$\{g_1,g_2,\ldots,g_{2r-1},g_{2r},
q_1,\ldots,q_{k-2r}\}=\{1,2,\ldots,k\}$ (see (\ref{leto5007after})),
$r=1,2,\ldots,$ $[k/2]$.

Assume that all traces
$\bar T^{k-2r}_{g_1,g_2,\ldots,g_{2r-1}, 
g_{2r}}\Phi(t_{q_1},\ldots,t_{q_{k-2r}})$
(for all $r=1,2,\ldots,$ $[k/2]$ and for all possible
$g_1,g_2,\ldots,g_{2r-1},g_{2r}$ (see {\rm (\ref{leto5007after})))
exist.

The analogue of the Hu--Meyer formula from \cite{SU11} for the multiple Stra\-to\-no\-vich 
stochastic integral (\ref{novem2026xxx10}) and the case of a multidimensional
Wiener process can be written as
$$
\bar J^{S}[\Phi]_{T,t}^{(i_1\ldots i_k)}
=
J'[\Phi]_{T,t}^{(i_1\ldots i_k)}+
$$

\vspace{-3mm}
$$
+
\sum\limits_{r=1}^{[k/2]}
\sum_{\stackrel{(\{\{g_1, g_2\}, \ldots, 
\{g_{2r-1}, g_{2r}\}\}, \{q_1, \ldots, q_{k-2r}\})}
{{}_{\{g_1, g_2, \ldots, 
g_{2r-1}, g_{2r}, q_1, \ldots, q_{k-2r}\}=\{1, 2, \ldots, k\}}}}
\prod\limits_{s=1}^r
{\bf 1}_{\{i_{g_{{}_{2s-1}}}=~i_{g_{{}_{2s}}}\ne 0\}}\times
$$

\vspace{4mm}
\begin{equation}
\label{novem2026xxx101}
\times
J'\hspace{-1mm}\left[\bar T^{k-2r}_{g_1,g_2,\ldots,g_{2r-1}, 
g_{2r}}\Phi\right]_{T,t}^{(i_{q_1}\ldots i_{q_{k-2r}})}\ \ \ \hbox{w.~p.~1,}
\end{equation}

\vspace{4mm}
\noindent
where $J'\hspace{-1mm}\left[\bar T^{k-2r}_{g_1,g_2,\ldots,g_{2r-1}, 
g_{2r}}\Phi\right]_{T,t}^{(i_{q_1}\ldots i_{q_{k-2r}})}\stackrel{\sf def}{=}
\bar T^{k-2r}_{g_1,g_2,\ldots,g_{2r-1},g_{2r}}\Phi$ for $k=2r,$ 
the existence of the right-hand side of (\ref{novem2026xxx101})
ensures the existence of the multiple Stratonovich stochastic
integral $\bar J^{S}[\Phi]_{T,t}^{(i_1\ldots i_k)};$
another notations are the same as in (\ref{2026may10}).

Assume that all limiting traces
$T^{k-2r}_{g_1,g_2,\ldots,g_{2r-1}, g_{2r}}\Phi (t_{q_1},\ldots,t_{q_{k-2r}})$
(see (\ref{novem2026xxx1122})) exist
for all $r=1,2,\ldots,$ $[k/2]$ and for all possible
$g_1,g_2,\ldots,g_{2r-1},g_{2r}$ (see {\rm (\ref{leto5007after})).
Moreover, suppose that
$$
\bar T^{k-2r}_{g_1,g_2,\ldots,g_{2r-1}, 
g_{2r}}\Phi(t_{q_1},\ldots,t_{q_{k-2r}})=
T^{k-2r}_{g_1,g_2,\ldots,g_{2r-1}, 
g_{2r}}\Phi(t_{q_1},\ldots,t_{q_{k-2r}})
$$

\newpage
\noindent
almost everywhere on $[t, T]^{k-2r}$ (with respect to Lebesgue's measure)
for all $r=1,2,\ldots,$ $[k/2]$ and for all possible
$g_1,g_2,\ldots,g_{2r-1},g_{2r}$ (see {\rm (\ref{leto5007after})),
where $\bar T^{k-2r}_{g_1,g_2,\ldots,g_{2r-1}, 
g_{2r}}\Phi(t_{q_1},\ldots,t_{q_{k-2r}})$
is defined by (\ref{novem2026xxx11}).

Then, using (\ref{2026may10}), (\ref{novem2026xxx101}), and (\ref{2026may2}), we obtain

\vspace{-3mm}
$$
{\sf M}\left\{\left(\bar J^{S}[\Phi]_{T,t}^{(i_1\ldots i_k)}-\hat J^{S}[\Phi]_{T,t}^{(i_1\ldots i_k)}  
\right)^2\right\}\le 
$$

$$
\le  
\bar C_k
\sum\limits_{r=1}^{[k/2]}
\sum_{\stackrel{(\{\{g_1, g_2\}, \ldots, 
\{g_{2r-1}, g_{2r}\}\}, \{q_1, \ldots, q_{k-2r}\})}
{{}_{\{g_1, g_2, \ldots, 
g_{2r-1}, g_{2r}, q_1, \ldots, q_{k-2r}\}=\{1, 2, \ldots, k\}}}}
\prod\limits_{s=1}^r
{\bf 1}_{\{i_{g_{{}_{2s-1}}}=~i_{g_{{}_{2s}}}\ne 0\}}\times
$$

$$
\Biggl.\times
\left\Vert \bar T^{k-2r}_{g_1,g_2,\ldots,g_{2r-1}, 
g_{2r}}\Phi-
T^{k-2r}_{g_1,g_2,\ldots,g_{2r-1}, 
g_{2r}}\Phi\right\Vert_{L_2([t, T]^{k-2r})}^2=0,
$$

\vspace{2mm}
\noindent
where $\bar C_k$ is a constant.

Then (see (\ref{novem2026xxx4}), (\ref{novem2026xxx5}))
$$
\bar J^{S}[\Phi]_{T,t}^{(i_1\ldots i_k)}=\hat J^{S}[\Phi]_{T,t}^{(i_1\ldots i_k)}=
\hbox{\vtop{\offinterlineskip\halign{
\hfil#\hfil\cr
{\rm l.i.m.}\cr
$\stackrel{}{{}_{p\to \infty}}$\cr
}} }
\sum_{j_1,\ldots,j_k=0}^{p}
C_{j_k\ldots j_1}
\prod_{l=1}^k \zeta_{j_l}^{(i_l)}
$$

\noindent
w.~p.~1, where 
$C_{j_k\ldots j_1}$ is the Fourier coefficient
defined by (\ref{2025may30}).

\subsection{Expansion of Iterated Stratonovich Stochastic Integrals of Arbitrary
Multiplicity}

{\bf Theorem~9.} {\it
Let $\{\phi_j(x)\}_{j=0}^{\infty}$ be an arbitrary 
CONS in $L^2([t, T]),$
and let $\psi_1(\tau) \equiv \dots \equiv \psi_k(\tau) \equiv 1$. 
For any multiplicity $k \in {\bf N}$ 
and any indices $i_1, \dots, i_k \in \{1, \dots, d\},$
the partial sum 
\begin{equation}
\label{ud11}
S_{T,t}^{(i_1, \dots, i_k)p}
=\sum_{j_1,\ldots,j_k=0}^{p}
C_{j_k\ldots j_1}
\prod_{l=1}^k \zeta_{j_l}^{(i_l)}
\end{equation}
with equal truncation limits $p_1 = \dots = p_k = p$ 
converges in the mean-square sense 
to the iterated Stratonovich stochastic integral 
$J^{*}[\psi^{(k)}]_{T,t}^{(i_1\ldots i_k)}$ as $p \to \infty:$

\newpage
\noindent
$$
\lim_{p \to \infty} 
{\sf M}\left\{
\left| S_{T,t}^{(i_1, \dots, i_k)p} 
- J^{*}[\psi^{(k)}]_{T,t}^{(i_1\ldots i_k)} \right|^2 
\right\} = 0,
$$
where $J^{*}[\psi^{(k)}]_{T,t}^{(i_1\ldots i_k)}$ is defined by {\rm (\ref{str})}
for $i_1,\ldots,i_k=1,\ldots,d,$ 
$C_{j_k \ldots j_1}$ has the form {\rm (\ref{after3000}),} 
and
$$
\zeta_{j}^{(i)}=
\int\limits_t^T \phi_{j}(\tau) d{\bf w}_{\tau}^{(i)}\ \ \ (i=1,\ldots,d,\ j=0,1,2,\ldots)
$$ 
are independent standard Gaussian random variables for various 
$i$ or $j.$ 
}

{\bf Proof.}\
Without loss of generality, 
by a linear time rescaling, 
we assume that $[t, T] = [0, 1]$, 
since all variation norms and estimates 
are invariant under linear time changes.

{\bf Step~1. Identification with the Karhunen--Lo\`eve lift.}
Let ${\bf w}_{\tau} = \bigl({\bf w}_{\tau}^{(1)}, \dots, {\bf w}_{\tau}^{(d)}\bigr)$, 
$\tau \in [0, 1]$, be a standard 
$d$-dimensional Wiener process with independent components and  ${\bf w}_0 = 0$ w.~p.~1. 
For each coordinate $\mathbf{w}^{(i)}$
($i = 1, \dots, d$), let $\mathcal{K}_i \subset L^2(\Omega)$ denote the
closed linear Gaussian subspace spanned by $\{\mathbf{w}_{\tau}^{(i)} : \tau \in [0,1]\}$.
Following Friz and Victoir (2010, Sect.~3.2) \cite{Friz}, the Cameron--Martin space
$\mathcal{H}_i$ associated with $\mathbf{w}^{(i)}$ is the Hilbert space of
paths $h : [0,1] \to \mathbf{R}$ of the form
$$
h(\tau) = \mathsf{M}\{Z \mathbf{w}_{\tau}^{(i)}\}, \quad \tau \in [0, 1]
$$
for some $Z \in \mathcal{K}_i$, equipped with the inner product
$$
\langle h, g \rangle_{\mathcal{H}_i} := \mathsf{M}\{Z Z'\},
$$
where $g(\tau) = \mathsf{M}\{Z' \mathbf{w}_{\tau}^{(i)}\}$ with $Z' \in \mathcal{K}_i$.

Since every element of $\mathcal{K}_i$ admits the Wiener integral representation
$$
Z = \int\limits_0^1 f_1(s)\, d\mathbf{w}_s^{(i)}, \ \ \ 
Z' = \int\limits_0^1 f_2(s)\, d\mathbf{w}_s^{(i)}
$$
for unique functions $f_1, f_2 \in L^2([0,1])$, the It\^o isometry yields
$$
h(\tau) = \mathsf{M}\left\{\int\limits_0^1 f_1(s)\, d\mathbf{w}_s^{(i)}
\int\limits_0^{\tau} 1\, d\mathbf{w}_s^{(i)}\right\} = \int\limits_0^{\tau} f_1(s)\, ds,
$$
and similarly 
$$
g(\tau) = \int\limits_0^{\tau} f_2(s) ds,
$$
with inner product
$$
\langle h, g \rangle_{\mathcal{H}_i}
= \mathsf{M}\left\{\int\limits_0^1 f_1(s)\, d\mathbf{w}_s^{(i)}
\int\limits_0^1 f_2(s)\, d\mathbf{w}_s^{(i)}\right\}
= \int\limits_0^1 f_1(s) f_2(s)ds.
$$

Defining the antiderivatives 
$$
h_j(\tau) = \int\limits_0^{\tau} \phi_j(s)ds,\ \ \
\tau \in [0,\,1],\ \ \  j \ge 0,
$$
we obtain 
$$
\langle h_j, h_m \rangle_{\mathcal{H}_i} 
= \int\limits_0^1 \phi_j(s)\phi_m(s)\,ds = \delta_{j,m}.
$$ 

Hence, $\{h_j\}_{j=0}^{\infty}$ forms 
an orthonormal basis of $\mathcal{H}_i$, 
and the vectors with $h_j$ in the $i$-th slot 
form an orthonormal basis of 
$$
\mathcal{H} = \bigoplus\limits_{i=1}^d \mathcal{H}_i.
$$

In Sect.~6 of Friz and Victoir (2010) \cite{Friz}, 
the general $L^2$-approximation of a Gaussian process $X_{\tau}$ 
for a subset $A \subset {\bf N}$ is defined by
$$
X_{\tau}^A := {\sf M}\bigl\{ X_{\tau} \mid \mathcal{F}_A \bigr\} 
= \sum_{j \in A} Z_j h_j(\tau),\ \ \
\mathcal{F}_A = \sigma(Z_j : j \in A).
$$

Setting $X_{\tau} = {\bf w}_{\tau}^{(i)}$, $A = A_p := \{0, 1, \dots, p\}$, 
and 
$$
Z_j = \zeta_j^{(i)} := \int\limits_0^1 \phi_j(s)d{\bf w}_{s}^{(i)},
$$ 
this specializes to the absolutely continuous 
componentwise approximation (also see (\ref{novem2026xxx1000}) for $t=0$)
$$
{\bf w}_{\tau}^{(i)p}
:= {\sf M}\bigl\{ {\bf w}_{\tau}^{(i)} \mid \mathcal{F}_{A_p}^{(i)} \bigr\}
= \sum_{j=0}^{p} \zeta_j^{(i)} h_j(\tau),\ \ \ 
i = 1, \dots, d,
$$
where $\mathcal{F}_{A_p}^{(i)} = \sigma\bigl(\zeta_j^{(i)} : 0 \le j \le p\bigr)$. 
We denote the resulting $d$-dimensional absolutely continuous approximation 
(having sample paths in $C^{1\text{-var}}([0,1],{\bf R}^d)$ w.~p.~1) by 
$$
{\bf w}_{\tau}^{p} := \bigl({\bf w}_{\tau}^{(1)p}, \dots, {\bf w}_{\tau}^{(d)p}\bigr),\ \ \
\tau \in [0, 1].
$$

Next, for any path $x \in C^{1\text{-var}}([0,\,1], {\bf R}^d)$, 
its canonical step-$N$ lift $S_N(x)$ 
into the truncated tensor algebra 
$$
T^N({\bf R}^d) = \bigoplus\limits_{m=0}^N ({\bf R}^d)^{\otimes m}
$$ 
is given by \cite{Friz} (Friz and Victoir, 2010, Sect.~1.1):
$$
S_N(x)_{\tau} 
= 1 + \sum_{m=1}^{N}\
\int\limits_{0 < s_1 < \dots < s_m < \tau} 
dx_{s_1} \otimes \dots \otimes dx_{s_m}.
$$

Let $\{e_1, \dots, e_d\}$ be the standard orthonormal basis of ${\bf R}^d$ 
given by $e_r = (\delta_{r,1}, \dots, \delta_{r,d})^{\top}$ ($r = 1, \dots, d$). 
Then, in the $k$-th tensor power 
$({\bf R}^d)^{\otimes k} = \underbrace{{\bf R}^d \otimes \dots \otimes {\bf R}^d}_{k \text{ times}}$ 
(of dimension $d^k$), 
the set $\{e_{j_1} \otimes \dots \otimes e_{j_k}\}_{j_1, \dots, j_k = 1}^d$ 
forms an orthonormal basis with respect to the canonical Euclidean structure, satisfying
$$
\bigl\langle e_{j_1} \otimes \dots \otimes e_{j_k},\, e_{i_1} \otimes \dots \otimes e_{i_k} \bigr\rangle 
= \prod_{m=1}^k \delta_{j_m, i_m}.
$$

Substituting $N = k$, $x = {\bf w}^p$, and $\tau = 1$ (noting that 
$S_k({\bf w}^p)_{0,1} := \left(S_k({\bf w}^p)_0\right)^{-1} \otimes S_k({\bf w}^p)_1 = S_k({\bf w}^p)_1$ 
since $S_k({\bf w}^p)_0 = 1$), the lift $S_k({\bf w}^p)_{0,1}$ in the direct sum 
$$
T^k({\bf R}^d) = \bigoplus\limits_{m=0}^k ({\bf R}^d)^{\otimes m}
$$ is written componentwise as
$$
\begin{aligned}
S_k({\bf w}^p)_{0,1} = 
\Biggl( 1,\; \int\limits_{0 < s_1 < 1} d{\bf w}^p_{s_1},\; 
&\int\limits_{0 < s_1 < s_2 < 1} d{\bf w}^p_{s_1} \otimes d{\bf w}^p_{s_2},\; \dots, \\
&\dots,\; \int\limits_{0 < s_1 < \dots < s_k < 1} d{\bf w}^p_{s_1} 
\otimes \dots \otimes d{\bf w}^p_{s_k} \Biggr).
\end{aligned}
$$

Applying the canonical projection $\pi_k : T^k({\bf R}^d) \to ({\bf R}^d)^{\otimes k}$ 
onto the $k$-th tensor level, defined explicitly 
for any tuple 
$$
a = (a_0, a_1, \dots, a_k) \in \bigoplus\limits_{m=0}^k ({\bf R}^d)^{\otimes m}
$$ 
by $\pi_k(a_0, a_1, \dots, a_k) = a_k$, we obtain the expansion
$$
\pi_k\bigl(S_k({\bf w}^p)_{0,1}\bigr) 
= \int\limits_{0 < s_1 < \dots < s_k < 1} 
d{\bf w}^p_{s_1} \otimes \dots \otimes d{\bf w}^p_{s_k}
= 
$$
$$
=\sum_{j_1, \dots, j_k = 1}^{d} 
S^{(j_1, \dots, j_k)p} \, e_{j_1} \otimes \dots \otimes e_{j_k},
$$
where $S^{(j_1, \dots, j_k)p}$ is defined by (\ref{ud11}) for $T=1,$ $t=0.$

Consequently, the partial sum $S^{(i_1, \dots, i_k)p}$ 
is the $(i_1, \dots, i_k)$-component of $\pi_k \bigl(S_k({\bf w}^p)_{0,1}\bigr)$:
$$
S^{(i_1, \dots, i_k)p} 
= \int\limits_{0 < s_1 < \dots < s_k < 1} 
d{\bf w}^{(i_1)p}_{s_1} \cdots d{\bf w}^{(i_k)p}_{s_k} 
= 
$$
$$
=\bigl\langle \pi_k\bigl(S_k({\bf w}^p)_{0,1}\bigr),\, 
e_{i_1} \otimes \dots \otimes e_{i_k} \bigr\rangle.
$$

\vspace{1mm}

\textbf{Step 2. Convergence at level~$2$ in $p'$-variation.} 
Fix a rough path exponent $p' \in (2, 3)$ (in the sense of Lyons, 1998, 
Sect.~1.2 and~2.2 \cite{Lyons}), so that $\lfloor p' \rfloor = 2$. 
Recall that for any $N \ge 1$ and geometric $p'$-rough paths 
$\mathbf{x}, \mathbf{y} \in C_0^{0,p'\text{-var}}\bigl(\lbrack 0, 1 \rbrack, G^N(\mathbf{R}^d)\bigr)$ 
\cite{Friz} (Friz and Victoir, 2010, Sect.~1.1, p.~372), the $p'$-variation 
distance is defined by
$$
d_{p'\text{-var};\lbrack 0, 1 \rbrack}(\mathbf{x}, \mathbf{y}) 
= \sup_{D=\{0=t_0<t_1<\dots<t_r=1\}} 
\left( \sum_{l=0}^{r-1} d(\mathbf{x}_{t_l, t_{l+1}}, \mathbf{y}_{t_l, t_{l+1}})^{p'} \right)^{1/p'},
$$
where $d(g, h) = \|g^{-1} \otimes h\|$ denotes the Carnot--Carath\'eodory 
metric on $G^N(\mathbf{R}^d)$ and $\mathbf{x}_{s,\tau} := \mathbf{x}_s^{-1} \otimes \mathbf{x}_\tau$.

For the $d$-dimensional Wiener process $\mathbf{w} = (\mathbf{w}^{(1)}, \dots, \mathbf{w}^{(d)})$, 
the covariance of each coordinate $R(s, \tau) = \min(s, \tau)$ has finite 2D $1$-variation 
($\rho = 1$, so $p' > 2\rho = 2$; see Friz and Victoir, 2010, Sect.~3.1.1 \cite{Friz}). 
By Theorem~35(i), Definition~36, and Theorem~51 (eq.~(19)) in Friz and Victoir (2010) \cite{Friz} 
applied to $X = \mathbf{w}$ and its projection approximations $X^{A_p} = \mathbf{w}^p$ 
($A_p = \{0, \dots, p\}$), projecting onto level~$2$ via the canonical contraction 
$\pi_{0,2} : G^3(\mathbf{R}^d) \to G^2(\mathbf{R}^d)$ (defined by 
$\pi_{0,2}(1, g_1, g_2, g_3) = (1, g_1, g_2)$) yields
$$
\lim_{p \to \infty} \mathsf{M}\Bigl\{ d_{p'\text{-var};\lbrack 0, 1 \rbrack}\bigl(S_2(\mathbf{w}^p), \mathbf{W}^{(2)}\bigr)^q \Bigr\} = 0 
\qquad \forall q \in \lbrack 1, \infty),
$$
where $\mathbf{W}^{(2)} := \pi_{0,2}(\mathbf{W}) \in C_0^{0,p'\text{-var}}\bigl(\lbrack 0, 1 \rbrack, G^2(\mathbf{R}^d)\bigr)$ 
and $\pi_{0,2}\bigl(S_3(\mathbf{w}^p)\bigr) = S_2(\mathbf{w}^p)$.

For any subinterval $\lbrack s, \tau \rbrack \subset \lbrack 0, 1 \rbrack$, choosing 
a partition $D$ containing $s$ and $\tau$ (e.g., $D = \{0, s, \tau, 1\}$) gives
$$
d\bigl(S_2(\mathbf{w}^p)_{s,\tau}, \mathbf{W}^{(2)}_{s,\tau}\bigr) 
\le d_{p'\text{-var};\lbrack 0, 1 \rbrack}\bigl(S_2(\mathbf{w}^p), \mathbf{W}^{(2)}\bigr) 
\xrightarrow[p \to \infty]{L^q(\Omega)} 0 
\qquad \forall q \in \lbrack 1, \infty).
$$

By the equivalence of homogeneous norms on $G^2(\mathbf{R}^d)$ \cite{Friz} (Friz and Victoir, 
2010, Sect.~1.1), 
it follows that for every $\lbrack s, \tau \rbrack \subset \lbrack 0, 1 \rbrack$, 
the second-level coordinates of $S_2(\mathbf{w}^p)_{s,\tau}$ converge in $L^2(\Omega)$ 
to those of $\mathbf{W}^{(2)}_{s,\tau}$. 

On the other hand, for any complete orthonormal 
system $\{\phi_j(x)\}_{j=0}^\infty$ in $L^2(\lbrack 0, 1 \rbrack)$ and every 
subinterval $[s, \tau]\subset [0, 1]$, the second-order 
iterated Riemann--Stieltjes integrals of the series approximations $\mathbf{w}^p$ 
(i.e. $S^{(i_1 i_2)p}_{\tau,s}$)
over $[s, \tau]$ converge in $L^2(\Omega)$ \cite{2018a} (Kuznetsov, 2026) to the $2$-fold iterated 
Stratonovich stochastic integrals over $[s, \tau]$ 
(it is necessary
to take $\psi_1(t_1)={\bf 1}_{\{s<t_1\}},$ $\psi_2(t_2)={\bf 1}_{\{s<t_2<\tau\}}$ in Theorem~2.44
\cite{2018a} 
and Theorem 1.16 ($k=2$) \cite{2018a} 
and use the relation (2.125) \cite{2018a} for matrix and integral traces). By uniqueness of the 
$L^2(\Omega)$-limit and path continuity, $\mathbf{W}^{(2)}$ coincides w.~p.~$1$ with 
the canonical Stratonovich Brownian rough path of step~$2$ on $\lbrack 0, 1 \rbrack$.

\vspace{1mm}

\textbf{Step 3. Extension to arbitrary multiplicity $k$.} 
By Step~2 and Markov's inequality, $d_{p'\text{-var};\lbrack 0, 1 \rbrack}\bigl(S_2(\mathbf{w}^p), \mathbf{W}^{(2)}\bigr) \xrightarrow{\mathsf{P}} 0$ 
as $p \to \infty$. Assume $k \ge 3$ (the case $k \le 2$ is already covered by Step~2). 
Since $p' \in (2, 3)$, we have $\lfloor p' \rfloor = 2$.

By Lyons' Extension Theorem \cite{Lyons} (Lyons, 1998, Theorem~2.2.1), 
the extension map
$$
\tilde{S}_k : C_0^{0,p'\text{-var}}\bigl(\lbrack 0, 1 \rbrack, G^2(\mathbf{R}^d)\bigr) 
\to C_0^{0,p'\text{-var}}\bigl(\lbrack 0, 1 \rbrack, G^k(\mathbf{R}^d)\bigr)
$$
is unique and continuous with respect to $d_{p'\text{-var};\lbrack 0, 1 \rbrack}$, 
and satisfies $\tilde{S}_k\bigl(S_2(\mathbf{w}^p)\bigr) = S_k(\mathbf{w}^p)$. 
Since the family of $k$-fold iterated Stratonovich stochastic integrals of $\mathbf{w}$
up to order $k$ has the same first two levels as $\mathbf{W}^{(2)}$, satisfies
Chen's identity, obeys the classical shuffle product rule (hence takes values in
$G^k(\mathbf{R}^d)$), and has finite $p'$-variation w.~p.~1, it is uniquely
determined on all higher levels $3, \dots, k$ by its second-level rough path
$\mathbf{W}^{(2)}$ (via Lyons' Extension Theorem
\cite{Lyons} (Lyons, 1998, Theorem~2.2.1)). Consequently,
$\pi_k(\tilde{S}_k(\mathbf{W}^{(2)})_{0,1}) \in (\mathbf{R}^d)^{\otimes k}$
coincides w.~p.~1 with the tensor $J^*$ of $k$-fold iterated Stratonovich
stochastic integrals of $\mathbf{w}$ over $[0,1]$:
$$
\pi_k\bigl(\tilde{S}_k(\mathbf{W}^{(2)})_{0,1}\bigr) 
= J^* 
= \sum_{j_1,\dots,j_k=1}^d J^{*(j_1,\dots,j_k)} e_{j_1} \otimes \dots \otimes e_{j_k},
$$
so that $J^{*(i_1,\dots,i_k)} = \bigl\langle \pi_k\bigl(\tilde{S}_k(\mathbf{W}^{(2)})_{0,1}\bigr), e_{i_1} \otimes \dots \otimes e_{i_k} \bigr\rangle$, 
where $J^*[\psi^{(k)}]_{T,t}^{(i_1 \dots i_k)}$ for $\psi_1 \equiv \dots \equiv \psi_k \equiv 1$ 
and $T = 1$, $t = 0$ is denoted by $J^{*(i_1,\dots,i_k)}$.

Applying the continuity of the extension map $\tilde{S}_k$ with respect to 
$d_{p'\text{-var};\lbrack 0, 1 \rbrack}$ (via the Continuous Mapping Theorem \cite{Mann}), we obtain
$$
d_{p'\text{-var};\lbrack 0, 1 \rbrack}\bigl(S_k(\mathbf{w}^p), \tilde{S}_k(\mathbf{W}^{(2)})\bigr) 
= d_{p'\text{-var};\lbrack 0, 1 \rbrack}
\Bigl(\tilde{S}_k\bigl(S_2(\mathbf{w}^p)\bigr), \tilde{S}_k\bigl(\mathbf{W}^{(2)}\bigr)\Bigr) 
\xrightarrow[p \to \infty]{\mathsf{P}} 0.
$$

To pass to individual tensor coordinates, recall from Sect.~1.1 of Friz and Victoir 
(2010) \cite{Friz} the equivalence of homogeneous norms on $G^k(\mathbf{R}^d)$ 
(where $\|g\| = d(1, g)$ is the Carnot--Carath\'eodory norm entering $d(g, h) = \|g^{-1} \otimes h\|$):
$$
\frac{1}{K_k} \max_{1 \le m \le k} \lvert\pi_m(g)\rvert^{1/m} 
\le \|g\| 
\le K_k \max_{1 \le m \le k} \lvert\pi_m(g)\rvert^{1/m}
$$
for some constant $K_k \ge 1$ and all $g \in G^k(\mathbf{R}^d)$. Choosing the coarse 
partition $D = \{0, 1\}$ in the definition of $d_{p'\text{-var};\lbrack 0, 1 \rbrack}$ for 
$\mathbf{x} = S_k(\mathbf{w}^p)$ and $\mathbf{y} = \tilde{S}_k(\mathbf{W}^{(2)})$ 
(since paths in $C_0^{0,p'\text{-var}}\bigl(\lbrack 0, 1 \rbrack, G^k(\mathbf{R}^d)\bigr)$ satisfy 
the initial condition $\mathbf{x}_0 = \mathbf{y}_0 = 1$, we have 
$\mathbf{x}_{0,1} = \mathbf{x}_0^{-1} \otimes \mathbf{x}_1 = 
\mathbf{x}_1$ and $\mathbf{y}_{0,1} = \mathbf{y}_1$), 
we obtain
$$
d(\mathbf{x}_{0,1}, \mathbf{y}_{0,1}) = 
\|\mathbf{x}_{0,1}^{-1} \otimes \mathbf{y}_{0,1}\| 
\le d_{p'\text{-var};\lbrack 0, 1 \rbrack}(\mathbf{x}, \mathbf{y}).
$$

Using the algebraic identity $\mathbf{y}_{0,1} - \mathbf{x}_{0,1} = 
\mathbf{x}_{0,1} \otimes (\mathbf{x}_{0,1}^{-1} \otimes \mathbf{y}_{0,1} - 1)$ 
and the multiplication rule
$$
\pi_k(\mathbf{a} \otimes \mathbf{b}) = \sum_{m=0}^k \pi_{k-m}(\mathbf{a}) \otimes \pi_m(\mathbf{b})
$$
in $T^k(\mathbf{R}^d)$ \cite{Lyons} (Lyons, 1998, eq.~(2.2)), we obtain
$$
\begin{aligned}
\pi_k(\mathbf{y}_{0,1}) - \pi_k(\mathbf{x}_{0,1}) 
&= \pi_k\Bigl(\mathbf{x}_{0,1} \otimes (\mathbf{x}_{0,1}^{-1} \otimes \mathbf{y}_{0,1} - 1)\Bigr) \\
&= \sum_{m=0}^k \pi_{k-m}(\mathbf{x}_{0,1}) \otimes \pi_m(\mathbf{x}_{0,1}^{-1} \otimes \mathbf{y}_{0,1} - 1) \\
&= \sum_{m=1}^k \pi_{k-m}(\mathbf{x}_{0,1}) \otimes \pi_m(\mathbf{x}_{0,1}^{-1} \otimes \mathbf{y}_{0,1}),
\end{aligned}
$$
since $\pi_0(\mathbf{x}_{0,1}^{-1} \otimes \mathbf{y}_{0,1} - 1) = 1 - 1 = 0$ and 
$\pi_m(g - 1) = \pi_m(g)$ for all $m \ge 1$.

Applying $\lvert\pi_m(g)\rvert \le K_k^m \|g\|^m$ (which follows from the lower bound in 
the norm equivalence above) to $g = \mathbf{x}_{0,1}^{-1} \otimes \mathbf{y}_{0,1}$ 
(so that $\|g\| = d(\mathbf{x}_{0,1}, \mathbf{y}_{0,1})$) and using the triangle inequality 
together with the submultiplicativity $\lvert u \otimes v \rvert \le \lvert u \rvert\,\lvert v \rvert$ of the tensor norm, 
we deduce the explicit estimate
$$
\bigl\lvert\pi_k(\mathbf{y}_{0,1}) - \pi_k(\mathbf{x}_{0,1})\bigr\rvert 
\le \sum_{m=1}^k K_k^m \lvert\pi_{k-m}(\mathbf{x}_{0,1})\rvert\, d(\mathbf{x}_{0,1}, \mathbf{y}_{0,1})^m 
\qquad \forall k \ge 1.
$$
Since $d_{p'\text{-var};\lbrack 0, 1 \rbrack}(\mathbf{x}, \mathbf{y}) = d_{p'\text{-var};\lbrack 0, 1 \rbrack}\bigl(S_k(\mathbf{w}^p), \tilde{S}_k(\mathbf{W}^{(2)})\bigr) \xrightarrow{\mathsf{P}} 0$, 
it follows that $d(\mathbf{x}_{0,1}, \mathbf{y}_{0,1}) \xrightarrow{\mathsf{P}} 0$ 
(and hence $d(\mathbf{x}_{0,1}, \mathbf{y}_{0,1})^m \xrightarrow{\mathsf{P}} 0$ for every $m = 1, \dots, k$). 
Since the estimate above holds at every level (with $\pi_0(\mathbf{x}_{0,1}) = 1$), 
it follows by induction on the tensor level that each lower-level projection 
$\pi_{k-m}(\mathbf{x}_{0,1})$ ($m \in \{1, \dots, k\}$) converges in probability 
to $\pi_{k-m}(\mathbf{y}_{0,1})$ and hence $\lvert\pi_{k-m}(\mathbf{x}_{0,1})\rvert$ is bounded 
in probability as $p \to \infty$. Since the product of a sequence bounded in 
probability and a sequence converging to $0$ in probability converges to $0$ 
in probability, each term in the sum and thus the entire right-hand side 
converges to $0$ in probability, so that $\lvert\pi_k(\mathbf{y}_{0,1}) - \pi_k(\mathbf{x}_{0,1})\rvert \xrightarrow{\mathsf{P}} 0$. 
Substituting back $\mathbf{x} = S_k(\mathbf{w}^p)$ and $\mathbf{y} = \tilde{S}_k(\mathbf{W}^{(2)})$, 
and noting that the $(i_1, \dots, i_k)$-coordinate of a tensor cannot exceed its 
Euclidean norm, we obtain for every multi-index $(i_1, \dots, i_k)$:
$$
\bigl\lvert S^{(i_1,\dots,i_k)p} - J^{*(i_1,\dots,i_k)} \bigr\rvert 
\le \bigl\lvert\pi_k(\mathbf{y}_{0,1}) - \pi_k(\mathbf{x}_{0,1})\bigr\rvert 
= 
$$
$$
=\bigl\lvert\pi_k\bigl(\tilde{S}_k(\mathbf{W}^{(2)})_{0,1}\bigr) - 
\pi_k\bigl(S_k(\mathbf{w}^p)_{0,1}\bigr)\bigr\rvert 
\xrightarrow[p \to \infty]{\mathsf{P}} 0.
$$

\medskip
\noindent\textbf{Step 4. Upgrading to mean-square convergence.}
By Proposition~21 in Friz and Victoir (2010) \cite{Friz}, 
the iterated integral $S^{(i_1, \dots, i_k)p}$ 
of the Gaussian process ${\bf w}^p$ (whose sample paths are absolutely continuous
and hence of finite $1$-variation w.~p.~1)
lies in the $k$-th (inhomogeneous) Wiener chaos $\mathcal{H}^{\le k}$. 
Since $\mathcal{H}^{\le k}$ is closed under convergence in probability, 
we also have $J^{*(i_1, \dots, i_k)} \in \mathcal{H}^{\le k}$, 
and thus the error 
$$
Z_p := S^{(i_1, \dots, i_k)p} - J^{*(i_1, \dots, i_k)} 
\in \mathcal{H}^{\le k}.
$$

Recall Lemma~22 in Friz and Victoir (2010) \cite{Friz}
(hypercontractivity on the inhomogeneous Wiener chaos $\mathcal{H}^{\le k}$): 
for any $Z \in \mathcal{H}^{\le k}$ and $1 < p_0 < q < \infty$,
$$ |Z|_{L^{p_0}(\Omega)} 
\le |Z|_{L^q(\Omega)} 
\le |Z|_{L^{p_0}(\Omega)} (k+1)(q-1)^{k/2} 
\max\bigl(1, (p_0 - 1)^{-k}\bigr).
$$

Furthermore, by Remark~24 in Friz and Victoir (2010) \cite{Friz}, 
on $\mathcal{H}^{\le k}$ convergence in probability 
is equivalent to convergence in $L^q(\Omega)$ 
for every $q \in (0, \infty)$:
$$
Z_p \xrightarrow[p \to \infty]{{\sf P}} 0 
\quad \Longleftrightarrow \quad 
\lim_{p \to \infty} |Z_p|_{L^q(\Omega)} = 0 
\qquad \forall\, q \in (0, \infty).
$$

Substituting $Z_p = S^{(i_1, \dots, i_k)p} - J^{*(i_1, \dots, i_k)}$ 
and $q = 2$ into this equivalence, 
and using Step~3 ($Z_p \xrightarrow{{\sf P}} 0$), 
we finally obtain
$$
\lim_{p \to \infty} |Z_p|_{L^2(\Omega)}^2 
= \lim_{p \to \infty} 
{\sf M}\left\{ 
\left| S^{(i_1, \dots, i_k)p} 
- J^{*(i_1, \dots, i_k)} \right|^2
\right\} = 0.
$$

Theorem~9 is proved. $\Box.$

\vspace{1mm}

{\bf Theorem~10\ (Extension to weighted iterated Stratonovich stochastic integrals).}
{\it Let $k \ge 1,$ 
let $\{\phi_j(x)\}_{j=0}^{\infty}$ be an arbitrary 
CONS in $L^2([0, 1]),$
$(i_1, \dots, i_k) \in \{1, \dots, d\}^k,$
and let $\psi_1(\tau), \dots,$ $\psi_k(\tau) \in C^1([0,\,1])$ be deterministic 
continuously differentiable weight functions. 
Define the weighted iterated stochastic integrals {\rm (}Riemann--Stiltjes and Stratonovich$,$ 
respectively{\rm )}
$$
S[\psi^{(k)}]^{(i_1, \dots, i_k)p}_{1,0}
:= \int\limits_{0 < t_1 < \dots < t_k < 1} 
\psi_1(t_1) \cdots \psi_k(t_k)\, 
d{\bf w}_{t_1}^{(i_1)p} \cdots d{\bf w}_{t_k}^{(i_k)p}, 
$$
$$
J^{*}[\psi^{(k)}]_{1,0}^{(i_1 \ldots i_k)}
:= \int\limits_{0 < t_1 < \dots < t_k < 1} 
\psi_1(t_1) \cdots \psi_k(t_k) 
\circ d{\bf w}_{t_1}^{(i_1)} \cdots \circ d{\bf w}_{t_k}^{(i_k)},
$$
where 
$$
S_{1,0}^{(i_1, \dots, i_k)p}
=\sum_{j_1,\ldots,j_k=0}^{p}
C_{j_k\ldots j_1}
\prod_{l=1}^k \zeta_{j_l}^{(i_l)};
$$

\noindent
another notations are the same as in Theorem~{\rm 9}. Then

\vspace{-3mm}
$$
\lim_{p \to \infty} 
{\sf M}\left\{
\left|S[\psi^{(k)}]^{(i_1, \dots, i_k)p}_{1,0}
- J^{*}[\psi^{(k)}]_{1,0}^{(i_1 \ldots i_k)} \right|^2 
\right\} = 0.
$$
}

{\bf Proof.} Recall that $\mathbf{w}^p = (\mathbf{w}^{(1)p}, \dots, \mathbf{w}^{(d)p})$ 
denotes the $p$-th partial sum expansion of the standard $d$-dimensional Wiener process 
$\mathbf{w}$ with independent components introduced in Theorem~9
(see (\ref{novem2026xxx1000}) for $t=0$). 
For each $i_m\in\{1,\ldots,d\},$ $m \in \{1, \dots, k\},$ and $t \in \lbrack 0, 1 \rbrack$, 
define the weighted process
$$
Y_t^{(i_m,m)} := \int\limits_0^t \psi_m(s)\, d\mathbf{w}_s^{(i_m)}
$$
and its absolutely continuous approximation
$$
Y_{t}^{(i_m,m)p} := \int\limits_0^t \psi_m(s)\, d\mathbf{w}_s^{(i_m)p},
$$
where $d{\bf w}_s^{(i)p}$ $(i=1,\ldots,d)$ is defined by (\ref{novem2026xxx1}).

Applying It\^o's formula to $F(t, {\bf w}_t^{(i_m)}) = \psi_m(t){\bf w}_t^{(i_m)}$ 
(where $F(t,x)=\psi_m(t)x,$ $\partial^2 F/\partial x^2 \equiv 0$) 
with $\mathbf{w}_0^{(i_m)} = 0$, we obtain w.~p.~1
$$
Y_t^{(i_m,m)} = \int\limits_0^t \psi_m(s)\, d\mathbf{w}_s^{(i_m)} 
= \psi_m(t)\mathbf{w}_t^{(i_m)} - \int\limits_0^t \psi_m'(s)\mathbf{w}_s^{(i_m)}\, ds,
$$
and similarly, by classical integration by parts for the absolutely continuous path
$\mathbf{w}_s^{(i_m)p}$ (of finite $1$-variation w.~p.~1),
$$
Y_{t}^{(i_m,m)p} = \psi_m(t)\mathbf{w}_t^{(i_m)p} - \int\limits_0^t \psi_m'(s)\mathbf{w}_s^{(i_m)p}\, ds\ \ \ 
\hbox{w.~p.~1}.
$$

Since $\psi_m(s)$ is deterministic and of bounded variation, the It\^o and Stratonovich 
differentials coincide 
($\psi_m(t)\,d\mathbf{w}_t^{(i_m)} = \psi_m(t) \circ d\mathbf{w}_t^{(i_m)} = \circ\, dY_t^{(i_m,m)}$), 
and $dY_{t}^{(i_m,m)p} = \psi_m(t)d{\bf w}_t^{(i_m,m)p}$. 
Hence the weighted iterated 
integrals rewrite as unweighted iterated integrals of the $k$-dimensional processes 
$Y^p_t = (Y^{(i_1,1)p}_t,$ $\ldots, Y^{(i_k,k)p}_t)$ and
$Y_t = (Y^{(i_1,1)}_t, \ldots, Y^{(i_k,k)}_t),$ $i_1,\ldots,i_k=1,\ldots,d$:
$$
S[\psi^{(k)}]^{(i_1, \dots, i_k)p}_{1,0}
= \int\limits_{0 < t_1 < \dots < t_k < 1} dY_{t_1}^{(i_1,1)p} \cdots dY_{t_k}^{(i_k,k)p}, 
$$
$$
J^{*}[\psi^{(k)}]_{1,0}^{(i_1 \ldots i_k)}
= \int\limits_{0 < t_1 < \dots < t_k < 1} \circ\, dY_{t_1}^{(i_1,1)} \cdots \circ dY_{t_k}^{(i_k,k)}.
$$

Note that the components $Y^{(i_m,m)}_t$ and $Y^{(i_l,l)}_t$ of
$Y_t = (Y^{(i_1,1)}_t,\ldots,Y^{(i_k,k)}_t)$
are independent whenever their first upper indices differ
($i_m \ne i_l$). Even though components with $i_m = i_l$
($m \ne l$) are correlated, the full covariance matrix
$R_Y$ has finite 2D $1$-variation ($\rho = 1$) on $[0,1]^2$ (see definition of 2D $1$-variation
in \cite{Friz}),
as shown below.

Now observe that $Y_t = (Y^{(i_1,1)}_t,\ldots,Y^{(i_k,k)}_t)$ is a centered $k$-dimensional 
Gaussian process whose covariance function $R_Y$ has finite 2D $1$-variation 
($\rho = 1$) on $[0, 1]^2$: indeed, by It\^o's isometry, for any 
rectangle $[s_1, s_2]\times [t_1, t_2] \subset [0, 1]^2$, 
the rectangular covariance increment of $(Y^{(i_m,m)}_t, Y^{(i_l,l)}_t)$ vanishes if $i_m \ne i_l$ 
and equals 
$$
\int\limits_{[s_1, s_2]\cap [t_1, t_2]} \psi_m(\tau)\psi_l(\tau) d\tau
$$ 
if $i_m = i_l$, so summing over any grid partition of $[0, 1]^2$ gives 
a total variation bounded by $\max\limits_{1 \le m \le k} \|\psi_m\|_\infty^2 < \infty$.

Next, let $\mathcal{F}_p := \sigma(\zeta_j^{(i)} : 0 \le j \le p,\, 1 \le i \le d)$ 
be the filtration generated by the first $p+1$ expansion coefficients of $\mathbf{w}_t$. 
Because $\{\phi_j(x)\}_{j=0}^\infty$ is a complete orthonormal basis of $L^2([0, 1])$, 
the coefficients $\zeta_j^{(i)}$ are independent standard Gaussian random variables 
and the $\sigma$-algebra 
$$
\mathcal{F}_\infty := \sigma\left(\bigcup\limits_{p=0}^\infty \mathcal{F}_p\right)
$$ 
coincides (up to ${\sf P}$-null sets)
with $\sigma(\mathbf{w}_t : t \in \lbrack 0, 1 \rbrack)$, so every $Y_t^{(i_m,m)}$ 
is $\mathcal{F}_\infty$-measurable. Moreover, since the tail series 
$$
\sum_{j=p+1}^\infty \zeta_j^{(i)} \int\limits_0^t \phi_j(s)ds
$$ 
is independent of $\mathcal{F}_p$ 
and has zero mean, we have w.~p.~1
$$
\mathbf{w}_t^{(i)p} = {\sf M}\left\{\mathbf{w}_t^{(i)} \mid \mathcal{F}_p\right\}
$$ 
for all $t \in \lbrack 0, 1 \rbrack$. By linearity of conditional expectation in the 
representations above, it follows that w.~p.~1
$$
Y_{t}^{(i_m,m)p} = {\sf M}\left\{Y_t^{(i_m,m)} \mid \mathcal{F}_p \right\}
$$
for all $m \in \{1, \dots, k\},$ $i_m\in\{1,\ldots,d\},$
and $t \in [0, 1]$.

In particular, $(Y_{t}^{(i_m,m)p})_{p \ge 0}$ is a closed centered Gaussian martingale with 
absolutely continuous sample paths (of finite $1$-variation w.~p.~1), 
and by L\'evy's Martingale Convergence Theorem, 
$Y_{t}^{(i_m,m)p} \to Y_t^{(i_m,m)}$ w.~p.~1 and in $L^2(\Omega)$ as $p \to \infty$ for every 
$t \in \lbrack 0, 1 \rbrack$.

Fix a rough path exponent $p' \in (2, 3)$ (so that $p' > 2\rho = 2$ and $\lfloor p' \rfloor = 2$). 
Applying Theorem~35(i), Definition~36, and Theorem~51 in Friz and Victoir (2010) \cite{Friz}
to the conditional expectation approximations $Y_{t}^p = {\sf M}\left\{Y_t \mid \mathcal{F}_p \right\}$ 
and projecting onto level~$2$ via the canonical 
contraction $\pi_{0,2} : G^3(\mathbf{R}^k) \to G^2(\mathbf{R}^k)$, 
we obtain
$$
\lim\limits_{p \to \infty} \mathsf{M}\Bigl\{ d_{p'\text{-var};\lbrack 0, 1 
\rbrack}\bigl(S_2(Y^p), \mathbf{Y}^{(2)}\bigr)^q \Bigr\} = 0 
\qquad \forall q \in \lbrack 1, \infty),
$$

\vspace{1mm}
\noindent
where $\mathbf{Y}^{(2)} := \pi_{0,2}(\mathbf{Y}) 
\in C_0^{0,p'\text{-var}}\bigl(\lbrack 0, 1 \rbrack, G^2(\mathbf{R}^k)\bigr)$. 

In particular, for every subinterval $\lbrack s, \tau \rbrack \subset \lbrack 0, 1 \rbrack$, 
the second-level coordinates of $S_2(Y^p)_{s,\tau}$ converge in $L^2(\Omega)$ to 
those of $\mathbf{Y}^{(2)}_{s,\tau}$. 

On the other hand, since the $(m, l)$-coordinate 
of $S_2(Y^p)_{s,\tau}$ is the weighted second-order iterated Riemann--Stieltjes integral 
$$
\int\limits_{s < u_1 < u_2 < \tau} 
\psi_m(u_1)\psi_l(u_2) d\mathbf{w}_{u_1}^{(i_m)p}\, d\mathbf{w}_{u_2}^{(i_l)p}
$$ 
with smooth weights $\psi_m(\tau), \psi_l(\tau) \in C^1(\lbrack 0, 1 \rbrack),$ then
\cite{2018a} (Kuznetsov, 2026)
guarantees its $L^2(\Omega)$-convergence 
to the $2$-fold iterated Stratonovich stochastic integral 
(it is necessary to take $\psi_1(t_1)=\tilde \psi_1(t_1){\bf 1}_{\{s<t_1\}},$ 
$\psi_2(t_2)=\tilde \psi_2(t_2){\bf 1}_{\{s<t_2<\tau\}}$
($\tilde \psi_m(\tau), \tilde \psi_l(\tau) \in C^1([0, 1]$) in Theorem~2.44 \cite{2018a}
and Theorem 1.16 ($k=2$) \cite{2018a}
and use the relation (2.125) \cite{2018a} for matrix and integral traces)
$$
\int\limits_{s < u_1 < u_2 < \tau} \psi_m(u_1)\psi_l(u_2) 
\circ d\mathbf{w}_{u_1}^{(i_m)} \circ d\mathbf{w}_{u_2}^{(i_l)} 
= \int\limits_{s < u_1 < u_2 < \tau} \circ\, dY_{u_1}^{(m)} \circ dY_{u_2}^{(l)}.
$$

By uniqueness of the $L^2(\Omega)$-limit and path continuity, $\mathbf{Y}^{(2)}$ 
coincides w.~p.~$1$ with the canonical Stratonovich rough path of $Y$ of step~$2$ 
on $\lbrack 0, 1 \rbrack$.

Finally, repeating Steps~3 and~4 of the proof of Theorem~9 verbatim for the unweighted 
iterated integrals of $Y^p_t$ and $Y_t$ in $\mathbf{R}^k$ with multi-index $(1, 2, \dots, k)$ 
(using the uniqueness and continuity of the Lyons extension $\tilde{S}_k$ in Lyons, 
1998, Theorem~2.2.1~\cite{Lyons}, and hypercontractivity on the inhomogeneous 
Wiener chaos $\mathcal{H}^{\le k}$ 
in Friz and Victoir, 2010, Lemma~22 and Remark~24 \cite{Friz})
yields
$$
\lim_{p \to \infty} 
{\sf M}\left\{
\left|S[\psi^{(k)}]^{(i_1, \dots, i_k)p}_{1,0}
- J^{*}[\psi^{(k)}]_{1,0}^{(i_1 \ldots i_k)} \right|^2 
\right\} = 0.
$$

Theorem~10 is proved. $\Box$

\vspace{2mm}

\textbf{Remark~1.}\ {\it Since the approximation errors
$$
Z_p := S_{T,t}^{(i_1,\dots,i_k)p} - J^*[\psi^{(k)}]_{T,t}^{(i_1\dots i_k)}
$$
in Theorem~{\rm 9} {\rm (}for $\psi_1(\tau) \equiv \dots \equiv \psi_k(\tau) \equiv 1${\rm )} and
$$
Z_p[\psi^{(k)}] := S[\psi^{(k)}]_{1,0}^{(i_1,\dots,i_k)p}
- J^*[\psi^{(k)}]_{1,0}^{(i_1\dots i_k)}
$$
in Theorem~{\rm 10} belong to the $k$-th inhomogeneous Wiener chaos
$\mathcal{H}^{\le k}$ and converge to zero in probability as $p \to \infty,$
Remark~{\rm 24} in Friz and Victoir {\rm (2010) \cite{Friz}} immediately implies that both
Theorems~{\rm 9} and~{\rm 10} remain valid for convergence in $L^q(\Omega)$ for every
$q \in (0, \infty)$ {\rm (}rather than only for $q = 2${\rm )}. Namely$,$ under the
assumptions of Theorems~{\rm 9} and~{\rm 10,} respectively$,$ we have
$$
\lim_{p\to\infty} \mathsf{M}\left\{\left|
S_{T,t}^{(i_1,\dots,i_k)p} - J^*[\psi^{(k)}]_{T,t}^{(i_1\dots i_k)}
\right|^q\right\} = 0 \qquad \forall q \in (0, \infty)
$$
and
$$
\lim_{p\to\infty} \mathsf{M}\left\{\left|
S[\psi^{(k)}]_{1,0}^{(i_1,\dots,i_k)p} - J^*[\psi^{(k)}]_{1,0}^{(i_1\dots i_k)}
\right|^q\right\} = 0 \qquad \forall q \in (0, \infty).
$$
}

{\bf Remark~2.}\ {\it Theorems {\rm 9} and {\rm 10} prove Hypotheses {\rm 2.1, 2.2,} and {\rm 2.5}
{\rm (}the case $p_1=\ldots=p_k=p$ and $\psi_1(\tau),\ldots,\psi_k(\tau)\in C^1([t, T])${\rm )}
from the monograph {\rm \cite{2018a}} {\rm (}Sect.~{\rm 2.5,} {\rm 2.28)}
on the mean-square convergence of expansions of iterated 
Stratonovich stochastic integrals of arbitrary multiplicity 
$(k \in {\bf N})$ with respect to components of a multidimensional 
Wiener process into generalized multiple Fourier series.
}

\vspace{15mm}

\linespread{0.7}

\end{document}